\documentclass{amsart}
\usepackage[T1]{fontenc}
\usepackage{lmodern}
\usepackage{textcomp}
\usepackage{amsmath}
\usepackage{amssymb}
\usepackage{mleftright}
\usepackage{xspace}
\usepackage[
left=2.6cm,
right=2.6cm,
top=3.0cm,
bottom=3.0cm
]{geometry}
\usepackage{microtype}
\usepackage{graphicx}
\usepackage{tikz}
\usetikzlibrary{arrows.meta}
\usepackage[table]{xcolor}
\usepackage{subcaption}
\usepackage{array}
\usepackage{tabularx}
\usepackage{booktabs}
\usepackage{multirow}
\usepackage{diagbox}
\usepackage{threeparttable}
\usepackage{adjustbox}
\usepackage[figuresright]{rotating}
\usepackage[ruled,vlined]{algorithm2e}
\usepackage[inline]{enumitem}
\usepackage[mathlines]{lineno}
\usepackage{placeins}
\theoremstyle{plain}
\newtheorem{assumption}{Assumption}[section]
\newtheorem{lemma}{Lemma}[section]
\newtheorem{theorem}{Theorem}[section]
\theoremstyle{definition}
\newtheorem{example}{Example}[section]
\theoremstyle{remark}

\numberwithin{figure}{section}
\numberwithin{table}{section}
\numberwithin{equation}{section}
\allowdisplaybreaks[4]
\usepackage[colorlinks]{hyperref}

\title[Preconditioned RDA method for Darcy flow]{A Preconditioned Reconstructed Discontinuous Approximation for Darcy Flow}

	\author[M. Cui]{Ming Cui}
  \address{School of Mathematics, Statistics and Mechanics, Beijing University of Technology, Beijing 100124, P.R. China}
	\email{mingcui@bjut.edu.cn}
	
	\author[Y.-X. Gao]{Yanxiang Gao}\address{School of Mathematics, Statistics and Mechanics, Beijing University of Technology, Beijing 100124, P.R. China}
	\email{yanxianggao@emails.bjut.edu.cn}
  
	\author[R. Li]{Ruo Li}
  \address{CAPT, LMAM and School of Mathematical Sciences, Peking University, Beijing 100871, P.R. China}
  \address{Chongqing Research Institute of Big Data, Peking University, Chongqing, 401121, P.R. China}
	\email{rli@math.pku.edu.cn}

	\author[Q.-C. Liu]{Qicheng Liu}
  \address{Hangzhou International Innovation Institute of Beihang University, Hangzhou 311115, P.R. China}
	\email{qcliu@buaa.edu.cn}
	
	\author[F.-Y. Yang]{Fanyi Yang}
  \address{School of Mathematics, Sichuan University, Chengdu 610065, P.R. China}
	\email{yangfanyi@scu.edu.cn}
	
	\author[S.-H. Zhao]{Shuhai Zhao}
  \address{Institute of Applied Physics and Computational Mathematics, Beijing 100094, P.R. China}
	\email{shuhai@pku.org.cn}

\newcommand{\avg}[1]{\{\!\!\{#1\}\!\!\}}
\newcommand{\jump}[1]{[\![#1]\!]}
\newcommand{\DGnorm}[1]{|\!|\!|#1|\!|\!|_{\mathrm{DG}}}
\newcommand{\Ltnorm}[2]{\|#1\|_{L^{2}(#2)}}
\allowdisplaybreaks[4]

\newcommand{\ds}{\,\mathrm{d}s}
\newcommand{\dx}{\,\mathrm{d}\boldsymbol{x}}

\newcommand{\Th}{\mathcal{T}_{h}}
\newcommand{\eh}{\mathcal{E}_{h}}
\newcommand{\ehi}{\mathcal{E}_{h}^{\mathrm{int}}}
\newcommand{\ehb}{\mathcal{E}_{h}^{\partial}}

\newcommand{\ud}{\boldsymbol{u}}
\newcommand{\pd}{p}
\newcommand{\vd}{\boldsymbol{v}}
\newcommand{\qd}{q}

\newcommand{\od}{\Omega}
\newcommand{\gd}{\Gamma}
\newcommand{\nd}{\boldsymbol{n}}

\newcommand{\uh}{\boldsymbol{u}_{h}}
\newcommand{\ph}{p_{h}}
\newcommand{\vh}{\boldsymbol{v}_{h}}
\newcommand{\qh}{q_{h}}

\newcommand{\Uh}{\boldsymbol{\mathrm{V}}_{h}^{m_{\boldsymbol{v}}}}
\newcommand{\Ph}{{\mathrm{Q}}_{h}^{m_{q}}}
\newcommand{\Uconst}{\mathbf{U}_{h}^{0}}
\newcommand{\Pconst}{\mathbb{U}_{h}^{0}}
\newcommand{\Pmean}{\mathbb{U}_{h,0}^{0}}

\newcommand{\Rmsca}{\mathcal{R}^{m}}
\newcommand{\Rmvec}{\boldsymbol{\mathcal{R}}^{m}}
\newcommand{\mv}{m_{\boldsymbol{v}}}
\newcommand{\mq}{m_{q}}
\newcommand{\Rv}{\boldsymbol{\mathcal{R}}^{m_{\boldsymbol{v}}}}
\newcommand{\Rq}{\mathcal{R}_{0}^{m_{q}}}
\newcommand{\Rqsca}{\mathcal{R}^{m_{q}}}

\newcommand{\MatPre}{\mathbf{A}_{0}}
\newcommand{\Mat}{\mathbf{A}}
\newcommand{\bfx}{\mathbf{x}}
\newcommand{\bfy}{\mathbf{y}}
\newcommand{\trans}{\mathrm{T}}

\newcommand{\x}{\boldsymbol{x}}
\newcommand{\uhi}{u_{h,i}}

\newcommand{\dimTh}{N_{e}}

\begin{document}
	\makeatletter
	\typeout{DEBUG textwidth=\the\textwidth}
	\typeout{DEBUG textheight=\the\textheight}
	\typeout{DEBUG paperwidth=\the\paperwidth}
	\typeout{DEBUG paperheight=\the\paperheight}
	\typeout{DEBUG rotdriver=\meaning\rot@LR}
	\makeatother
	% \linenumbers*
	\raggedbottom
	% vim:spell:tw=70:fo+=Mn:cc=70
\begin{abstract}
	In this paper, we develop an efficient numerical method for Darcy flow based 
	on the reconstructed discontinuous approximation (RDA). 
	The key contribution is a mesh-robust piecewise-constant preconditioner for arbitrarily high-order 
	RDA discretizations of Darcy flow. 
	The method discretizes the Darcy equations by a stabilized discontinuous Galerkin scheme 
	in coupled velocity-pressure form, with RDA spaces for both velocity and pressure. 
	Function spaces of this type can achieve arbitrarily high-order approximation accuracy while retaining 
	only one degree of freedom per element. We establish an optimal-order error estimate
	in the energy norm and a suboptimal error estimate in the $L^{2}$ norm. 
	Benefiting from the property that the number of degrees of freedom of the RDA space is independent of the approximation order, 
	we can construct a simple and efficient preconditioner based on piecewise constant functions for algebraic systems 
	arising from arbitrarily high-order discretizations. 
	We further prove a mesh-independent spectral equivalence between the preconditioner and the original 
	algebraic system, so that the convergence of the GMRES residual does not slow down under mesh refinement.
	Numerical experiments in two and three dimensions confirm the theoretical convergence rates 
	and demonstrate the effectiveness and efficiency of the proposed preconditioning method.

  \noindent \textbf{keywords}:
		Darcy flow;
		Reconstructed discontinuous approximation;
		Discontinuous Galerkin method;
		Preconditioning.

\end{abstract}

	\maketitle
	% vim:spell:tw=70:fo+=Mn:cc=70
\section{Introduction}\label{sec:Introduction}

We consider the following Darcy flow problem in porous media, governed by Darcy's law and the mass conservation equation:
\begin{subequations}
	\label{eq:Darcy-Model}
	\begin{alignat}{2}
		\label{eq:momentum}
		\ud+\mathbf{K}\nabla\pd&=0,&\quad&\mathrm{in}\;\od,\\
		\label{eq:mass-balance}
		\nabla\cdot\ud&={f},&\quad&\mathrm{in}\;\od,\\
		\label{eq:B.C.}
		\ud\cdot\nd&=g,&\quad&\mathrm{on}\;\gd,
	\end{alignat}
\end{subequations}
where $\ud$ is the velocity, $\pd$ is the pressure, and $\mathbf{K}=\kappa\mathbf{I}_{d}$ denotes the isotropic permeability tensor, with $\mathbf{I}_{d}$ being the $d$-dimensional identity matrix. We assume that $\kappa\in{L}^{\infty}(\od)$, $\kappa^{-1}\in{W}^{1,\infty}(\od)$, and that there exist constants $\kappa_{\min}$ and $\kappa_{\max}$ such that $0<\kappa_{\min}\leq\kappa\leq\kappa_{\max}<\infty$. The functions $f\in{L}^{2}(\od)$ and $g\in{L}^{2}(\gd)$ denote the given source/sink term and boundary condition, respectively, and satisfy the compatibility condition
\begin{equation}\label{eq:compatibility-condition}
	\int_{\od}f\dx=\int_{\gd}g\ds.
\end{equation}
Since the pressure is determined only up to an additive constant, we impose the normalization condition
\begin{equation}\label{eq:p-uniqueness}
	\int_{\od}\pd\dx=0.
\end{equation}
This model has broad applications in groundwater seepage, petroleum recovery, geothermal energy, carbon dioxide sequestration, and flow simulations in fractured reservoirs~\cite{Christie2001SPE10,Steefel2022Darcy,Traverso2013MixRT}. Therefore, the development of efficient and robust numerical methods for Darcy flow problems has received considerable attention.

Finite element methods are widely used for the numerical solution of Darcy flow. 
According to how the velocity and pressure are coupled, they fall into two classes: 
the primal pressure form and the mixed form.
In the primal pressure form, the pressure is obtained by solving a second-order elliptic problem, 
and the velocity is then computed from the pressure gradient through Darcy's law. 
Although this formulation is easy to implement, the resulting velocity approximation may suffer 
a loss of accuracy and generally fails to satisfy local mass 
conservation~\cite{Cockburn2007CGFluxDarcy,Masud2002StabilizedDarcy}. 
Standard mixed finite element methods solve the velocity and pressure as a coupled system~\cite{Brezzi2006MixedDarcy,Traverso2013MixRT}. 
However, the coupled system is typically indefinite, so the velocity and pressure spaces must satisfy 
an inf-sup stability condition~\cite{BrezziFortin1991Mixed}. 
This severely restricts the choice of the discrete approximation spaces for velocity and pressure. 
By introducing additional stabilization terms, these restrictions can be relaxed to some extent. 
We refer the reader to~\cite{Badia2010Darcy,Bochev2006PPPDarcy,Correa2008UnconditionallyDarcy,Masud2002StabilizedDarcy,Nakshatrala2006StableDarcy} 
for some stabilized continuous Galerkin formulations and 
to~\cite{Badia2010Darcy,Brezzi2005MixDGDarcy,Hughes2006StabilizedDarcy} 
for some stabilized discontinuous Galerkin (DG) methods. 
In particular, DG methods offer several advantages, including local conservation properties, 
flexibility in handling complex geometries and heterogeneous coefficients, and ease of 
implementing high-order approximations. These methods originated in the numerical solution 
of hyperbolic conservation laws~\cite{CockburnShu2001RKDG,ReedHill1973DG}. 
They have since been widely applied to a variety of problems, such as 
elliptic problems~\cite{ArnoldBrezziCockburnMarini2002DG,SuliMozolevski2007Biharmonic} 
and saddle-point problems~\cite{Brezzi2005MixDGDarcy,ChenFan2026poroelasticity,CockburnKanschatSchotzauSchwab2002LDGStokes,
Hughes2006StabilizedDarcy}. 
However, DG discretizations often involve 
a large number of degrees of freedom, giving rise to large, ill-conditioned linear systems that are 
computationally expensive to solve.

To address these challenges, various solvers and preconditioners have been developed, 
including iterative preconditioning schemes~\cite{Allen1992Iterative}, Schur-complement 
solvers~\cite{Maryska2000Schur}, block-diagonal preconditioners~\cite{Baerland2020PrecondDarcy}, 
block-triangular preconditioners~\cite{Axelsson2015PrecondDarcy}, 
block-factorization preconditioners~\cite{Nardean2021PrecondDarcy}, 
and multigrid methods~\cite{Brenner2018MGDarcy}. Depending on the resulting system, 
Krylov subspace methods include the conjugate gradient method for symmetric positive definite 
systems~\cite{Cliffe2000ParallelMFEM,Dougherty1990PCG,Traverso2013MixRT}, 
the minimum residual method for symmetric indefinite s
ystems~\cite{Maryska1996MINRES,Pazner2024SaddleHdiv,Traverso2013MixRT}, and 
the generalized minimal residual (GMRES) method for nonsymmetric systems~\cite{Nardean2021PrecondDarcy}. 
Although preconditioning can improve the spectral properties and accelerate convergence, 
the large size of the discrete systems still entails substantial computational and storage costs.

In recent years, the reconstructed discontinuous approximation (RDA) method has emerged as a promising technique for constructing high-order approximation spaces with minimal degrees of freedom. The method was originally developed for elliptic problems~\cite{LiRuo2019discontinuous} and later extended to elliptic interface problems~\cite{LiRuo2026PrecondInterface,LiRuo2020interface}. The key idea is to construct a high-order approximation space by solving local constrained least squares fitting problems on element patches, resulting in an approximation space that possesses only one degree of freedom per element. Compared with standard DG methods, RDA can achieve comparable accuracy with substantially fewer degrees of freedom, particularly at higher approximation orders. Moreover, this construction naturally enables the establishment of a norm equivalence between the high-order space and the lowest-order space. Based on this equivalence property, optimal preconditioners can be constructed from the piecewise constant space, yielding condition numbers of the corresponding preconditioned systems that are independent of the mesh size. Numerical experiments further demonstrate mesh-robust iterative convergence and the computational efficiency of these preconditioners~\cite{LiRuo2024Preconditioned,LiRuo2026PrecondInterface}.

The main contributions of this paper are as follows. First, we employ the RDA framework to construct high-order velocity and pressure approximation spaces. The reconstructed spaces support approximations of arbitrarily high order while retaining only one degree of freedom per element for each scalar component. Second, inspired by~\cite{Correa2008UnconditionallyDarcy,Hughes2006StabilizedDarcy}, we develop a stabilized DG scheme without requiring the inf-sup condition and establish an optimal-order error estimate in the energy norm and a suboptimal-order error estimate in the $L^{2}$ norm. Third, we propose an efficient preconditioner for the coupled velocity--pressure system by extending the RDA preconditioning strategy originally introduced in the context of scalar elliptic problems~\cite{LiRuo2024Preconditioned}. The preconditioner is induced by the energy norm on the piecewise constant space, and its block-diagonal structure permits separate solves for the velocity and pressure blocks. Mesh-independent field-of-values bounds for the preconditioned system are established through the equivalence between the energy norms on the high- and low-order spaces. Numerical experiments in two and three dimensions show that preconditioned GMRES iteration counts remain nearly constant under mesh refinement at a fixed polynomial degree and that approximate block solvers further reduce solve times.

The remainder of this paper is organized as follows. In Section~2, we introduce the notation and the reconstructed discontinuous
approximation operators. In Section~3, we present the stabilized DG formulation for the Darcy problem and establish the corresponding error analysis. In Section~4, we construct a preconditioner for the discrete system and provide the theoretical justification. In Section~5, two- and three-dimensional numerical experiments are presented to demonstrate the effectiveness of the proposed numerical method and preconditioner. Finally, we conclude the paper in Section~6.

	% vim:spell:tw=70:fo+=Mn:cc=70
\section{Preliminaries}\label{sec:Recons}
In this section, we introduce the notation and the reconstruction operators used throughout the paper.

Let $\od\subset\mathbb{R}^{d}\,(d=2,3)$ be a bounded convex polygonal or polyhedral domain with boundary $\gd=\partial\od$. In this paper, for any subdomain $D\subset\od$, we denote by $H^{s}(D)$ the standard Sobolev space equipped with the norm $\|\cdot\|_{H^{s}(D)}$ and seminorm $|\cdot|_{H^{s}(D)}$. We define the zero-mean subspace of $L^2(\od)$ by $L_{0}^{2}(\od):=\{\qd\in{L}^{2}(\od)\;|\;\int_{\od}\qd\dx=0\}$. For vector-valued functions in $[H^{s}(D)]^{d}$, the corresponding norms and seminorms are defined componentwise. For a sufficiently smooth vector field $\vd=(v_{1},\ldots,v_{d})^{\trans}$, the curl operator is defined by
\begin{align*}
	\nabla\times\vd=\frac{\partial{v}_{2}}{\partial{x}}-\frac{\partial{v}_{1}}{\partial{y}},\qquad{d}=2,
\end{align*}
and
\begin{align*}
	\nabla\times\vd=\begin{pmatrix}
		\frac{\partial{v}_{3}}{\partial{y}}-\frac{\partial{v}_{2}}{\partial{z}}\\
		\frac{\partial{v}_{1}}{\partial{z}}-\frac{\partial{v}_{3}}{\partial{x}}\\
		\frac{\partial{v}_{2}}{\partial{x}}-\frac{\partial{v}_{1}}{\partial{y}}
	\end{pmatrix},\qquad{d}=3.
\end{align*}

Let $\Th$ be a shape-regular, quasi-uniform polygonal or polyhedral partition of $\od$ consisting of $\dimTh$ elements. Denote the sets of all interior and boundary $(d-1)$-dimensional faces by $\ehi$ and $\ehb$, respectively, and set $\eh=\ehi\cup\ehb$. For any $K\in\Th$ and $e\in\eh$, we denote by $h_{K}$ and $h_{e}$ the diameters of $K$ and $e$, respectively, and set $h=\max_{K\in\Th}h_{K}$.

Throughout this paper, $C$ and $c$ denote generic positive constants independent of the mesh size $h$ but possibly dependent on other model and discretization parameters.

Under the above assumptions, the following trace and inverse inequalities hold for scalar-valued functions and, componentwise, for vector-valued functions. For any $K\in\Th$, we have
\begin{align}
	\label{ineq:trace}
	\Ltnorm{\qd}{\partial{K}}^{2}&\leq{c}\big(h_{K}^{-1}\Ltnorm{\qd}{K}^{2}+h_{K}\|{\qd}\|_{H^{1}(K)}^{2}\big),\quad\forall\,\qd\in{H}^{1}(K), \\
	\label{ineq:inverse}
	\|{\qd}\|_{H^{r}(K)}&\leq{c}h_{K}^{s-r}\|\qd\|_{H^{s}(K)},\quad0\leq{s}\leq{r},\;\forall\,\qd\in\mathbb{P}_{m}(K),
\end{align}
where $\mathbb{P}_{m}(K)$ denotes the space of polynomials of degree no greater than $m$ on $K$.

In the framework of discontinuous Galerkin methods, we introduce the following average and jump operators for scalar-valued and vector-valued functions. For any interior face $e\in\ehi$ shared by two neighboring elements $K^{+}$ and $K^{-}$, let $\nd^{+}$ and $\nd^{-}$ denote the unit outward normal vectors on $\partial{K}^{+}$ and $\partial{K}^{-}$, respectively. We define
\begin{alignat*}{2}
	\avg{\qd}&:=\frac{1}{2}\big(\qd|_{\partial{K}^{+}}+\qd|_{\partial{K}^{-}}\big),&\quad\jump{\qd}&:=\qd|_{\partial{K}^{+}}\nd^{+}+\qd|_{\partial{K}^{-}}\nd^{-}, \\
	\avg{\vd}&:=\frac{1}{2}\big(\vd|_{\partial{K}^{+}}+\vd|_{\partial{K}^{-}}\big),&\quad\jump{\vd\otimes\nd}&:=\vd|_{\partial{K}^{+}}\otimes\nd^{+}+\vd|_{\partial{K}^{-}}\otimes\nd^{-},
	\\\jump{\vd\cdot\nd}&:=\vd|_{\partial{K}^{+}}\cdot\nd^{+}+\vd|_{\partial{K}^{-}}\cdot\nd^{-}.
\end{alignat*}
For any boundary face $e\in\ehb$ with the unit outward normal vector $\nd$, we define
\begin{alignat*}{2}
	\avg{\qd}&:=\qd,&\quad\jump{\qd}&:=\qd\nd, \\
	\avg{\vd}&:=\vd,
	&\quad\jump{\vd\cdot\nd}&:=\vd\cdot\nd.
\end{alignat*}

Based on the above definitions, we introduce the broken Sobolev space
\begin{align*}
	H^{s}(\od,\Th):=\left\{\qd\in{L}^{2}(\od)\,|\,\qd|_{K}\in{H}^{s}(K),\forall\,K\in\Th\right\}.
\end{align*}
The same convention is used for vector-valued functions in $[H^{s}(\od,\Th)]^{d}$.

In this paper, we adopt the reconstruction operators introduced in~\cite{LiRuo2019Biharmonic,LiRuo2018Stokes}. For completeness, we recall their definitions as follows.

For each element $K\in\Th$, we recursively construct an element patch $S(K)$ with a prescribed cardinality $\#S(K)$. Starting from $S_{0}(K)=\{K\}$, we successively enlarge the patch by adding all elements neighboring the current patch. More precisely, for $\ell=1,2,\ldots$, we define
\begin{align*}
	S_{\ell}(K)=\left\{K''\in\Th\,|\,\overline{K''}\cap\overline{K'}\neq\emptyset\text{ for some }K'\in S_{\ell-1}(K)\right\}.
\end{align*}
Once $\#S_{\ell}(K)\geq\#S(K)$, we select the $\#S(K)$ elements in $S_{\ell}(K)$ closest to $K$ and denote the resulting set by $S(K)$.

We denote by $\x_{K}$ the barycenter of the element $K$, and define $I(K)=\{\x_{K'}\;|\;K'\in{S}(K)\}$. For a given function $\qd\in\Pconst$, where $\Pconst$ denotes the scalar-valued piecewise constant space, the local reconstruction operator $\Rmsca_{K}$ is defined by solving the following least squares problem over the polynomial space $\mathbb{P}_{m}(S(K))$:
\begin{align*}
	\Rmsca_{K}\qd&:=\mathop{\arg\min}\limits_{p\in\mathbb{P}_{m}(S(K))}\sum_{\x\in{I}(K)}|p(\x)-\qd(\x)|^{2},\\
	&\text{subject to }p(\x_{K})=\qd(\x_{K}).
\end{align*}
To ensure the uniqueness of the solution to the above least squares problem, we impose the following unisolvence assumption on the collocation points in $I(K)$.
\begin{assumption}
	\label{assum:recons-uni}
	For any element $K\in\Th$ and polynomial $p\in\mathbb{P}_{m}(S(K))$, $p|_{I(K)}=0$ implies $p|_{S(K)}=0$.
\end{assumption}
Under Assumption~\ref{assum:recons-uni}, the local reconstruction operator $\Rmsca_{K}$ is well defined. Subsequently, the global reconstruction operator $\Rmsca$ is defined elementwise by
\begin{equation}\label{eq:recon}
	(\Rmsca\qd)|_{K}=(\Rmsca_{K}\qd)|_{K},\quad\forall\,K\in\Th.
\end{equation}

Let $\{e_{K}\}_{K\in\Th}$ be the basis of the piecewise constant space $\Pconst$, where
\begin{align*}
	e_{K}(\x)=
	\begin{cases}
		1,&\x\in K,\\0,&\x\not\in K.
	\end{cases}
\end{align*}
The corresponding reconstructed basis functions are defined by $\psi_{K}=\Rmsca{e}_{K},\,K\in\Th$. $\{\psi_{K}\}_{K\in\Th}$ forms a basis of the reconstructed space $\Rmsca\Pconst$. Therefore, the number of degrees of freedom of the reconstructed space is equal to the number of elements in $\Th$. For any $\qd\in\Pconst$, its reconstruction can be written as
\begin{align*}
	\Rmsca\qd=\sum_{K\in\Th}\qd(\x_{K})\psi_{K}.
\end{align*}
This representation naturally extends the reconstruction operator $\Rmsca$ to $C^{0}(\od)$. For vector-valued functions, the reconstruction operator is defined componentwise, and we use the notation $\Rmvec$.

For notational convenience, we introduce the following space, which will be used in the subsequent analysis. For a given integer $\mv>0$, we define the vector-valued reconstructed space
\begin{align*}
	\Uh:=\Rv\Uconst,
\end{align*}
where $\Uconst$ is the vector-valued piecewise constant space. To facilitate the subsequent analysis, we also introduce a zero-mean version of the reconstructed space. Let $\Pi_{0}:\,L^{2}(\od)\rightarrow{L}^{2}_{0}(\od)$ be defined by
\begin{align*}
	\Pi_{0}{\qd}=\qd-\frac{1}{|\od|}\int_{\od}\qd\dx.
\end{align*}
It is easy to verify that $\Pi_{0}$ is the $L^{2}$-orthogonal projection onto ${L}^{2}_{0}(\od)$. The mean-corrected reconstruction operator is defined by $\Rmsca_{0}:=\Pi_{0}\Rmsca$. For a given integer $\mq>0$, we introduce the zero-mean piecewise constant space and the zero-mean reconstructed space as follows:
\begin{align*}
	\Pmean:=\Pconst\cap{L}_{0}^{2}(\od),\quad\Ph:=\Rq\Pmean\subset{L}^{2}_{0}(\od).
\end{align*}
Thus, both $\Pmean$ and $\Ph$ have dimension $\dimTh-1$. Accordingly, functions in these spaces admit coefficient representations in $\mathbb{R}^{\dimTh-1}$.

For the above scalar-valued and vector-valued reconstruction operators, the corresponding approximation properties are given as follows~\cite{LiRuo2018Stokes}.
\begin{lemma}\label{lem:Recons-error-estimate-scalar}
	Let $t\ge1$ be an integer and $\qd\in{H}^{t+1}(\od)\cap{L}^{2}_{0}(\od)$. Then there exists a constant $C$, independent of $h$, such that
	\begin{align}\label{ineq:recons-scalar}
		\|\qd-\Rq{\qd}\|_{H^{s}(K)}&\leq{C}\overline{\Lambda}_{\mq}h^{r_{\mq}+1-s}\left(\|\qd\|_{H^{r_{\mq}+1}(S(K))}+h^{s+d/2}\|\qd\|_{H^{r_{\mq}+1}(\od)}\right),\\
		\label{ineq:recons-scalar-trace}
		\|\nabla^{s}(\qd-\Rq{\qd})\|_{L^{2}(\partial{K})}&\leq {C}\overline{\Lambda}_{\mq}h^{r_{\mq}-s+1/2}\left(\|\qd\|_{H^{r_{\mq}+1}(S(K))}+h^{s+d/2}\|\qd\|_{H^{r_{\mq}+1}(\od)}\right),
	\end{align}
	where $0\leq{s}\leq{}r_{\mq}$, $r_{\mq}=\min\{t,{\mq}\}$, and the scalar reconstruction constant is defined by
	\begin{align*}
		\overline{\Lambda}_{\mq}
		=\max_{K\in\Th}\bigg\{1+\sqrt{\#S(K)}
		\bigg(\max\limits_{0\neq{\pd}\in\mathbb{P}_{\mq}(S(K))}\frac{\|\pd\|_{L^{2}(K)}^2}{h_{K}^{d}\sum\limits_{K'\in{S}(K)}|\pd(\x_{K'})|^{2}}\bigg)^{1/2}\bigg\}.
	\end{align*}
\end{lemma}
Compared with the corresponding estimates for $\Rq$ in~\cite{LiRuo2018Stokes}, the estimates for $\Rq$ contain an additional global-norm term introduced by the zero-mean projection $\Pi_{0}$. However, this additional term does not change the convergence orders in the subsequent analysis.
\begin{lemma}\label{lem:Recons-error-estimate-vector}
	Let $t\ge1$ be an integer and $\vd\in[{H}^{t+1}(\od)]^{d}$. Then there exists a constant $C$, independent of $h$, such that
	\begin{align}\label{ineq:recons-vector}
		\|\vd-\Rv\vd\|_{H^{s}(K)}&\leq{C}\overline{\Lambda}_{\mv}h^{r_{\mv}+1-s}\|\vd\|_{H^{r_{\mv}+1}(S(K))},\\
		\label{ineq:recons-vector-trace}
		\|\nabla^{s}(\vd-\Rv\vd)\|_{L^{2}(\partial{K})}&\leq{C}\overline{\Lambda}_{\mv}h^{r_{\mv}-s+1/2}\|\vd\|_{H^{r_{\mv}+1}(S(K))},
	\end{align}
	where $0\leq{s}\leq{}r_{\mv}$, $r_{\mv}=\min\{t,{\mv}\}$, and the constant $\overline{\Lambda}_{\mv}$ is defined as in Lemma~\ref{lem:Recons-error-estimate-scalar}, with $\mq$ replaced by $\mv$.
\end{lemma}
Following \cite{LiRuo2024Preconditioned,LiRuo2026PrecondInterface,LiRuo2019discontinuous}, we assume that the element patches are sufficiently large to ensure that the reconstruction constants $\overline{\Lambda}_{\mq}$ and $\overline{\Lambda}_{\mv}$ are uniformly bounded with respect to $h$. In this paper, for each fixed $m\in\{\mq,\mv\}$, we choose the patch cardinality such that $\#S(K)\geq\dim(\mathbb{P}_{m}(K))+m+1$ for $d=2$ and $\#S(K)\geq1.5\dim(\mathbb{P}_{m}(K))$ for $d=3$.

	% vim:spell:tw=70:fo+=Mn:cc=70
\section{Numerical Scheme}\label{sec:Darcy}

In this section, we develop a stabilized discontinuous Galerkin discretization for the Darcy problem and establish its error estimates.

We first introduce several bilinear forms used in the formulation. Let $\alpha,\beta>0$ be penalty parameters and set $\xi:=\frac{1}{\kappa}\in{W}^{1,\infty}(\od)$. For any $\ud,\vd\in[{H}^{1}(\od,\Th)]^{d},\,\pd,\qd\in{H}^{1}(\od,\Th)\cap{L}^{2}_{0}(\od)$, we define
\begin{flalign*}
	a_{h}(\ud,\vd)&=\sum\limits_{K\in\Th}\int_{K}\xi\ud\cdot\vd\dx,&&\\
	s_{\ud,h}(\ud,\vd)&=\frac{1}{2}\sum\limits_{K\in\Th}\int_{K}\xi(\nabla\cdot\ud)(\nabla\cdot\vd)\dx+\frac{1}{2}\sum\limits_{K\in\Th}\int_{K}\kappa(\nabla\times(\xi\ud))\cdot(\nabla\times(\xi\vd))\dx&&\\&\quad+\sum\limits_{e\in\ehi}\int_{e}\frac{\alpha}{h_{e}}\jump{\xi^{\frac{1}{2}}\ud\otimes\nd}:\jump{\xi^{\frac{1}{2}}\vd\otimes\nd}\ds+\sum\limits_{e\in\ehb}\int_{e}\frac{\alpha}{h_{e}}\jump{\xi^{\frac{1}{2}}\ud\cdot\nd}\jump{\xi^{\frac{1}{2}}\vd\cdot\nd}\ds,&&\\
	b_{h}(\vd,\pd)&=-\sum\limits_{K\in\Th}\int_{K}\pd\nabla\cdot\vd\dx+\sum\limits_{e\in\eh}\int_{e}\avg{\pd}\jump{\vd\cdot\nd}\ds,&&\\
	s_{\pd,h}(\pd,\qd)&=\sum\limits_{e\in\ehi}\int_{e}\frac{\beta}{h_{e}}\jump{\kappa^{\frac{1}{2}}\pd}\cdot\jump{\kappa^{\frac{1}{2}}\qd}\ds.&&
\end{flalign*}
For the second term in $s_{\ud,h}$, the displayed expression applies when $d=3$. For $d=2$, this term takes the form $\frac{1}{2}\sum_{K\in\Th}\int_{K}\kappa(\nabla\times(\xi\ud))(\nabla\times(\xi\vd))\dx$. This convention is used throughout the paper. Motivated by~\cite{Correa2008UnconditionallyDarcy,Hughes2006StabilizedDarcy,Masud2002StabilizedDarcy}, we employ a stabilized mixed formulation as the underlying discretization. The formulation contains consistent residual-based stabilization associated with Darcy's law, the mass conservation equation, and the curl of Darcy's law, together with face penalty terms imposed on both interior and boundary faces.
Accordingly, for all $\ud,\vd\in[{H}^{1}(\od,\Th)]^{d}$ and $\pd,\qd\in{H}^{1}(\od,\Th)\cap{L}^{2}_{0}(\od)$, we define the bilinear form by
\begin{equation}\label{eq:bilinear-formulation}
	\mathcal{A}(\ud,\pd;\vd,\qd)=a_{h}(\ud,\vd)+s_{\ud,h}(\ud,\vd)+b_{h}(\vd,\pd)
	-b_{h}(\ud,\qd)-\frac{1}{2}\sum\limits_{K\in\Th}\int_{K}\kappa(\xi\ud+\nabla\pd)\cdot(\xi\vd-\nabla\qd)\dx+s_{\pd,h}(\pd,\qd)
\end{equation}
and the linear form is defined by
\begin{equation}\label{eq:linear-formulation}
	\mathcal{F}(\vd,\qd)=\frac{1}{2}\sum\limits_{K\in\Th}\int_{K}\xi{f}\nabla\cdot\vd\dx+\sum\limits_{e\in\ehb}\int_{e}\frac{\alpha\xi}{h_{e}}g\nd\cdot\vd\ds
	+\sum\limits_{K\in\Th}\int_{K}f\qd\dx-\sum\limits_{e\in\ehb}\int_{e}{g}\qd\ds.
\end{equation}
Then the corresponding discrete problem for system~\eqref{eq:Darcy-Model} is as follows: find $(\uh,\ph)\in\Uh\times\Ph$ such that
\begin{equation}\label{eq:Darcy-discrete-problem}
	\mathcal{A}(\uh,\ph;\vh,\qh)=\mathcal{F}(\vh,\qh),\quad\forall\,(\vh,\qh)\in\Uh\times\Ph.
\end{equation}
\begin{lemma}\label{lem:consistency}
	Let $(\ud,\pd)\in[{H}^{1}(\od)]^{d}\times\left({H}^{1}(\od)\cap{L}^{2}_{0}(\od)\right)$ be the exact solution of the Darcy equations~\eqref{eq:Darcy-Model}. Then $(\ud,\pd)$ satisfies
	\begin{equation}\label{eq:consistency}
		\mathcal{A}(\ud,\pd;\vh,\qh)=\mathcal{F}(\vh,\qh),\quad\forall\,(\vh,\qh)\in\Uh\times\Ph.
	\end{equation}
\end{lemma}
By this lemma, subtracting~\eqref{eq:Darcy-discrete-problem} from~\eqref{eq:consistency} yields the error equation
\begin{equation}\label{eq:error-equation}
	\mathcal{A}(\ud-\uh,\pd-\ph;\vh,\qh)=0,\quad\forall\,(\vh,\qh)\in\Uh\times\Ph.
\end{equation}

For the subsequent analysis, we first present the following equivalence estimate.
\begin{lemma}
	\label{lem:norm-equivalence-velocity}
	Let $\alpha>0$ be a constant, and assume that $\xi\in{W}^{1,\infty}(\od)$ satisfies $0<\xi_{\min}\leq\xi\leq\xi_{\max}<\infty$. Then there exist constants $c,C>0$, depending on $\kappa$, $\xi$, and $\alpha$, but independent of $h$, such that for all $\uh\in\Uh$, the following inequality holds:
	\begin{align*}
		&\quad{c}\bigg[\sum\limits_{K\in\Th}\left(\Ltnorm{\xi^{\frac{1}{2}}\uh}{K}^{2}+\Ltnorm{\xi^{\frac{1}{2}}\nabla\cdot\uh}{K}^{2}+\Ltnorm{\kappa^{\frac{1}{2}}\nabla\times(\xi\uh)}{K}^{2}\right)\\&\qquad+\sum\limits_{e\in\ehi}\frac{\alpha}{h_{e}}\Ltnorm{\jump{\xi^{\frac{1}{2}}\uh\otimes\nd}}{e}^{2}+\sum\limits_{e\in\ehb}\frac{\alpha}{h_{e}}\Ltnorm{\jump{\xi^{\frac{1}{2}}\uh\cdot\nd}}{e}^{2}\bigg]
		\\&\leq\sum\limits_{K\in\Th}\left(\Ltnorm{\xi^{\frac{1}{2}}\uh}{K}^{2}+\Ltnorm{\xi^{\frac{1}{2}}\nabla\uh}{K}^{2}\right)+\sum\limits_{e\in\ehi}\frac{\alpha}{h_{e}}\Ltnorm{\jump{\xi^{\frac{1}{2}}\uh\otimes\nd}}{e}^{2}+\sum\limits_{e\in\ehb}\frac{\alpha}{h_{e}}\Ltnorm{\jump{\xi^{\frac{1}{2}}\uh\cdot\nd}}{e}^{2}\\
		&\leq{C}\bigg[\sum\limits_{K\in\Th}\left(\Ltnorm{\xi^{\frac{1}{2}}\uh}{K}^{2}+\Ltnorm{\xi^{\frac{1}{2}}\nabla\cdot\uh}{K}^{2}+\Ltnorm{\kappa^{\frac{1}{2}}\nabla\times(\xi\uh)}{K}^{2}\right)\\&\qquad+\sum\limits_{e\in\ehi}\frac{\alpha}{h_{e}}\Ltnorm{\jump{\xi^{\frac{1}{2}}\uh\otimes\nd}}{e}^{2}+\sum\limits_{e\in\ehb}\frac{\alpha}{h_{e}}\Ltnorm{\jump{\xi^{\frac{1}{2}}\uh\cdot\nd}}{e}^{2}\bigg].
	\end{align*}
\end{lemma}
\begin{proof}
	The lower bound follows directly from the assumptions on $\kappa$ and $\xi$, together with the relations among the divergence, curl, and gradient. Therefore, it remains only to prove the upper bound. Moreover, since the same terms appear on both sides, it suffices to show that $\sum_{K\in\Th}\Ltnorm{\xi^{\frac{1}{2}}\nabla\uh}{K}^{2}$ can be controlled by the terms on the right-hand side. In view of the assumption on $\xi$ and the identity $\nabla\times(\xi\uh)=\xi\nabla\times\uh+\nabla\xi\times\uh$, it suffices to prove
	\begin{equation}
		\begin{aligned}
			\label{ineq:norm-eq-H1}
			\sum\limits_{K\in\Th}\Ltnorm{\nabla\uh}{K}^{2}&\leq{C}\bigg(\sum\limits_{K\in\Th}\big(\Ltnorm{\nabla\cdot\uh}{K}^{2}+\Ltnorm{\nabla\times\uh}{K}^{2}\big)\\&\qquad+\sum\limits_{e\in\ehi}\frac{1}{h_{e}}\Ltnorm{\jump{\uh\otimes\nd}}{e}^{2}+\sum\limits_{e\in\ehb}\frac{1}{h_{e}}\Ltnorm{\jump{\uh\cdot\nd}}{e}^{2}\bigg).
		\end{aligned}
	\end{equation}
	
	Let $\tilde{\Th}$ be a conforming, shape-regular simplicial refinement of $\Th$, obtained by subdividing each polygonal or polyhedral element into triangles or tetrahedra, respectively. We consider the $\mathcal{C}^{0}$ Lagrange finite element space $\mathbf{S}_{h}^{\mv}=[\mathbb{P}_{\mv}(\tilde{\Th})]^{d}\cap[{H}^{1}(\od)]^{d}$. To describe functions in this space, for each element $\tilde{K}\in\tilde{\Th}$, we select a set of Lagrange nodes $\mathcal{N}_{\tilde{K}}=\{\x_{\tilde{K},i}\}_{i=1}^{\mathrm{dim}_{\mv}}$, where $\mathrm{dim}_{\mv}=\dim(\mathbb{P}_{\mv}(\tilde{K}))$, and denote the set of all Lagrange nodes by $\mathcal{N}=\bigcup_{\tilde{K}\in\tilde{\Th}}\mathcal{N}_{\tilde{K}}$. Let $\phi_{\tilde{K},i}$ denote the Lagrange basis functions associated with the nodes $\x_{\tilde{K},i}$. Then every $\vh\in\mathbf{S}_{h}^{\mv}$ can be expressed as
	\begin{align*}
		\vh|_{\tilde{K}}=\sum_{i=1}^{\mathrm{dim}_{\mv}}\vh(\x_{\tilde{K},i})\phi_{\tilde{K},i},\quad{\tilde{K}}\in\tilde{\Th}.
	\end{align*}
	Similarly, the discrete function $\uh$ admits the following local representation on each element:
	\begin{align*}
		\uh|_{\tilde{K}}=\sum_{i=1}^{\mathrm{dim}_{\mv}}\uh|_{\tilde{K}}(\x_{\tilde{K},i})\phi_{\tilde{K},i},\quad{\tilde{K}}\in\tilde{\Th}.
	\end{align*}
	Following~\cite[Lemma~4]{LiRuo2021SequentialLeastSquares}, we first classify the nodes according to their locations.
	\begin{align*}
		\mathcal{N}^{\mathrm{int}}&=\{\x\in\mathcal{N}\;|\;\x\text{ is in the interior of an element }K\in\Th\},\\
		\mathcal{N}^{\mathrm{s}}&=\{\x\in\mathcal{N}\;|\;\x\text{ is common to several elements of }\Th,\,\x\notin\gd\},\\
		\mathcal{N}^{\mathrm{b}}&=\{\x\in\mathcal{N}\;|\;\x\in\gd\}.
	\end{align*}
	For each node $\x$, we define the sets
	\begin{align*}
		\mathcal{T}_{\x}&=\{K\in\Th\;|\;\x\in\overline{K}\},\\
		\mathcal{T}^{\partial}_{\x}&=\{K\in\Th\;|\;\exists{e}\in\ehb\text{ such that }e\subset\partial{K},\,\x\in\overline{e}\},\\\mathcal{V}^{\partial}_{\x}& =\operatorname*{span}\{\nd\;|\;{e}\in\ehb,\,\x\in\overline{e},\,\nd\text{ is the unit outward normal vector on }e\}.
	\end{align*}
	For each boundary node $\x\in\mathcal{N}^{\mathrm{b}}$, let
	$\{\nd_{\x,i}\}_{i=1}^{\dim\mathcal{V}_{\x}^\partial}$ be an orthonormal basis of $\mathcal{V}^{\partial}_{\x}$. For each interior face $e\in\ehi$, let $K_{e}^{+}$ and $K_{e}^{-}$ denote the two elements sharing $e$. For each boundary face $e\in\ehb$, let $K_{e}$ denote the element such that $e\subset\partial{K}_{e}$, and let $\nd_{e}$ denote the unit outward normal vector on $e$. Based on these sets, we define the continuous finite element function $\vh$ by specifying its values as follows:
	\begin{align*}
		\vh(\x)=
		\begin{cases}
			\uh(\x),\,&\x\in\mathcal{N}^{\mathrm{int}},\\
			\frac{\sum_{K\in\mathcal{T}_{\x}}\uh|_{K}(\x)}{\#\mathcal{T}_{\x}},\,&\x\in\mathcal{N}^{\mathrm{s}},\\
			\frac{\sum_{K\in\mathcal{T}^{\partial}_{\x}}\uh|_{K}(\x)}{\#\mathcal{T}^{\partial}_{\x}}-\sum\limits_{i=1}^{\dim\mathcal{V}_{\x}^\partial}\left(\frac{\sum_{K\in\mathcal{T}^{\partial}_{\x}}\uh|_{K}(\x)}{\#\mathcal{T}^{\partial}_{\x}}\cdot\nd_{\x,i}\right)\nd_{\x,i},\,&\x\in\mathcal{N}^{\mathrm{b}}.
		\end{cases}
	\end{align*}
	
	Therefore, we decompose $\uh$ as $\uh=(\uh-\vh)+\vh$. By the triangle inequality, we obtain
	\begin{align*}
		\sum\limits_{K\in\Th}\Ltnorm{\nabla\uh}{K}^{2}\leq{C}\bigg(\sum\limits_{K\in\Th}\Ltnorm{\nabla(\uh-\vh)}{K}^{2}+\Ltnorm{\nabla\vh}{\od}^{2}\bigg).
	\end{align*}
	By the definitions of $\mathcal{V}_{\x}^\partial$ and $\vh$, for every boundary face $e\in\ehb$, $\vh(\x)\cdot\nd_{e}=0$ at all Lagrange nodes $\x\in\overline{e}$, and hence $\vh\cdot\nd_{e}=0$ on $e$. Therefore, $\vh$ satisfies the homogeneous normal boundary condition $\vh\cdot\nd=0$ on $\gd$, and thus $\vh\in{H}_{0}(\operatorname{div};\od)\cap{H}(\operatorname{curl};\od)$. Applying the inequalities~\cite[Theorem~3.1]{CaiZhiqiang2000LSStokes} for $d=2$ and~\cite[Chapter~I, Theorem~3.9]{Girault1986FEMNS} for $d=3$, we obtain
	\begin{equation}
		\label{ineq:Gaffney}
		\Ltnorm{\nabla\vh}{\od}^{2}\leq{C}\left(\Ltnorm{\nabla\cdot\vh}{\od}^{2}+\Ltnorm{\nabla\times\vh}{\od}^{2}\right).
	\end{equation}
	
	Next, we estimate $\Ltnorm{\nabla(\uh-\vh)}{K}$. For any $K\in\Th$, we have
	\begin{align*}
		\Ltnorm{\nabla(\uh-\vh)}{K}^{2}&\leq{C}\sum_{i=1}^{\mathrm{dim}_{\mv}}\sum\limits_{\tilde{K}\in\tilde{\Th},\tilde{K}\subset{K}}\Ltnorm{\nabla\phi_{\tilde{K},i}}{\tilde{K}}^{2}|\uh(\x_{\tilde{K},i})-\vh(\x_{\tilde{K},i})|^{2}\\
		&\leq{C}h_{K}^{d-2}\bigg(\sum\limits_{\x\in\mathcal{N}^{\mathrm{s}}\cap\overline{K}}\big|(\uh-\vh)|_{K}(\x)\big|^{2}+\sum\limits_{\x\in\mathcal{N}^{\mathrm{b}}\cap\overline{K}}\big|(\uh-\vh)|_{K}(\x)\big|^{2}\bigg).
	\end{align*}

	For any $\x\in\mathcal{N}^{\mathrm{s}}\cap\overline{K}$,
	\begin{align*}
		|\uh(\x)-\vh(\x)|^{2}&=\bigg|\frac{\sum_{K'\in\mathcal{T}_{\x}}\uh|_{K'}(\x)}{\#\mathcal{T}_{\x}}-\frac{\#\mathcal{T}_{\x}\uh|_{K}(\x)}{\#\mathcal{T}_{\x}}\bigg|^{2}\leq{C}\sum\limits_{e\in\ehi,\x\in\overline{e}}\big|\uh|_{K_{e}^{+}}(\x)-\uh|_{K_{e}^{-}}(\x)\big|^{2}\\&\leq{C}\sum\limits_{e\in\ehi,\x\in\overline{e}}\|\jump{\uh\otimes\nd}\|_{L^{\infty}(e)}^{2}.
	\end{align*}
	For any $\x\in\mathcal{N}^{\mathrm{b}}\cap\overline{K}$,
	\begin{align*}
		|\uh(\x)-\vh(\x)|^{2}&=\bigg|\frac{\sum_{K'\in\mathcal{T}^{\partial}_{\x}}\uh|_{K'}(\x)}{\#\mathcal{T}^{\partial}_{\x}}-\frac{\#\mathcal{T}^{\partial}_{\x}\uh|_{K}(\x)}{\#\mathcal{T}^{\partial}_{\x}}-\sum\limits_{i=1}^{\dim\mathcal{V}_{\x}^\partial}\left(\frac{\sum_{K'\in\mathcal{T}^{\partial}_{\x}}\uh|_{K'}(\x)}{\#\mathcal{T}^{\partial}_{\x}}\cdot\nd_{\x,i}\right)\nd_{\x,i}\bigg|^{2}\\
		&\leq{C}\bigg(\sum\limits_{e\in\ehi,\x\in\overline{e}}\left|\uh|_{K_{e}^{+}}(\x)-\uh|_{K_{e}^{-}}(\x)\right|^{2}+\sum\limits_{e\in\ehb,\x\in\overline{e}}\left|\uh|_{K_{e}}(\x)\cdot\nd_{e}\right|^{2}\bigg)\\
		&\leq{C}\bigg(\sum\limits_{e\in\ehi,\x\in\overline{e}}\|\jump{\uh\otimes\nd}\|_{L^{\infty}(e)}^{2}+\sum\limits_{e\in\ehb,\x\in\overline{e}}\|\jump{\uh\cdot\nd}\|_{L^{\infty}(e)}^{2}\bigg).
	\end{align*}
	Summing the estimates over all elements, we have
	\begin{equation}
		\begin{aligned}
			\label{ineq:nabla-u-v}
			\sum\limits_{K\in\Th}\Ltnorm{\nabla(\uh-\vh)}{K}^{2}&\leq{C}\bigg(\sum\limits_{e\in\ehi}h_{e}^{d-2}\|\jump{\uh\otimes\nd}\|_{L^{\infty}(e)}^{2}+\sum\limits_{e\in\ehb}h_{e}^{d-2}\|\jump{\uh\cdot\nd}\|_{L^{\infty}(e)}^{2}\bigg)\\
			&\leq{C}\bigg(\sum\limits_{e\in\ehi}h_{e}^{-1}\|\jump{\uh\otimes\nd}\|_{L^{2}(e)}^{2}+\sum\limits_{e\in\ehb}h_{e}^{-1}\|\jump{\uh\cdot\nd}\|_{L^{2}(e)}^{2}\bigg).
		\end{aligned}
	\end{equation}
	By the triangle inequality,
	\begin{align*}
		\Ltnorm{\nabla\cdot\vh}{\od}^{2}+\Ltnorm{\nabla\times\vh}{\od}^{2}\leq{C}\sum\limits_{K\in\Th}\left(\Ltnorm{\nabla\cdot\uh}{K}^{2}+\Ltnorm{\nabla\times\uh}{K}^{2}+\Ltnorm{\nabla(\uh-\vh)}{K}^{2}\right).
	\end{align*}
	Combining~\eqref{ineq:Gaffney} and~\eqref{ineq:nabla-u-v},~\eqref{ineq:norm-eq-H1} follows. The proof is complete.
\end{proof}
Based on this estimate, for any $\uh\in[{H}^{1}(\od,\Th)]^{d}$ and $\ph\in{H}^{1}(\od,\Th)\cap{L}^{2}_{0}(\od)$, we define the DG energy norms $\DGnorm{\cdot}$ by
\begin{align*}
	\DGnorm{\uh}=\bigg(&\sum\limits_{K\in\Th}\big(\Ltnorm{\xi^{\frac{1}{2}}\uh}{K}^{2}+\Ltnorm{\xi^{\frac{1}{2}}\nabla\uh}{K}^{2}\big)\\&+\sum\limits_{e\in\ehi}\frac{\alpha}{h_{e}}\Ltnorm{\jump{\xi^{\frac{1}{2}}\uh\otimes\nd}}{e}^{2}+\sum\limits_{e\in\ehb}\frac{\alpha}{h_{e}}\Ltnorm{\jump{\xi^{\frac{1}{2}}\uh\cdot\nd}}{e}^{2}\bigg)^{\frac{1}{2}},\\
	\DGnorm{\ph}=\bigg(&\sum\limits_{K\in\Th}\Ltnorm{\kappa^{\frac{1}{2}}\nabla\ph}{K}^{2}+\sum\limits_{e\in\ehi}\frac{\beta}{h_{e}}\Ltnorm{\jump{\kappa^{\frac{1}{2}}\ph}}{e}^{2}\bigg)^{\frac{1}{2}}.
\end{align*}
With the above definitions, for any $(\uh,\ph)\in[{H}^{1}(\od,\Th)]^{d}\times\left({H}^{1}(\od,\Th)\cap{L}^{2}_{0}(\od)\right)$, we define the energy norm $\DGnorm{(\cdot,\cdot)}$ on the product space by
\begin{equation}\label{eq:DG-norm}
	\begin{aligned}
		\DGnorm{(\uh,\ph)}=\big(\DGnorm{\uh}^{2}+\DGnorm{\ph}^{2}\big)^{\frac{1}{2}}.
	\end{aligned}
\end{equation}
Then we establish the boundedness and coercivity of the bilinear form with respect to this norm. Combining Lemma~\ref{lem:norm-equivalence-velocity} with integration by parts, the trace inequality, and the Cauchy--Schwarz inequality, we obtain the following boundedness result.
\begin{lemma}\label{lem:continuity}
	Let $\alpha,\beta>0$. Then there exists a constant $C>0$ independent of $h$ such that the following estimate holds for any $(\uh,\ph)\in[{H}^{1}(\od,\Th)]^{d}\times\left({H}^{1}(\od,\Th)\cap{L}^{2}_{0}(\od)\right)$ and $(\vh,\qh)\in\Uh\times\Ph$:
	\begin{equation}\label{ineq:continuity}
		|\mathcal{A}(\uh,\ph;\vh,\qh)|\leq{C}\DGnorm{(\uh,\ph)}\DGnorm{(\vh,\qh)}.
	\end{equation}
\end{lemma}
\begin{lemma}\label{lem:coercive}
	Assume that the penalty parameters $\alpha$ and $\beta$ are positive. Then there exists a positive constant $C$, independent of $h$, such that for any $(\uh,\ph)\in\Uh\times\Ph$:
	\begin{equation}\label{ineq:coercive}
		\mathcal{A}(\uh,\ph;\uh,\ph)\geq{C}\DGnorm{(\uh,\ph)}^{2}.
	\end{equation}
\end{lemma}
\begin{proof}
	For any $(\uh,\ph)\in\Uh\times\Ph$, by the definition of the bilinear form~\eqref{eq:bilinear-formulation} and Lemma~\ref{lem:norm-equivalence-velocity}, we obtain
	\begin{align*}
		&\quad\mathcal{A}(\uh,\ph;\uh,\ph)
		\\&=\frac{1}{2}\sum\limits_{K\in\Th}\int_{K}\xi\uh\cdot\uh\dx+\frac{1}{2}\sum\limits_{K\in\Th}\int_{K}\xi(\nabla\cdot\uh)(\nabla\cdot\uh)\dx+\frac{1}{2}\sum\limits_{K\in\Th}\int_{K}\kappa(\nabla\times(\xi\uh))\cdot(\nabla\times(\xi\uh))\dx\\&\quad+\sum\limits_{e\in\ehi}\int_{e}\frac{\alpha}{h_{e}}\jump{\xi^{\frac{1}{2}}\uh\otimes\nd}:\jump{\xi^{\frac{1}{2}}\uh\otimes\nd}\ds+\sum\limits_{e\in\ehb}\int_{e}\frac{\alpha}{h_{e}}\jump{\xi^{\frac{1}{2}}\uh\cdot\nd}\jump{\xi^{\frac{1}{2}}\uh\cdot\nd}\ds\\&\quad+\frac{1}{2}\sum\limits_{K\in\Th}\int_{K}\kappa\nabla\ph\cdot\nabla\ph\dx+\sum\limits_{e\in\ehi}\int_{e}\frac{\beta}{h_{e}}\jump{\kappa^{\frac{1}{2}}\ph}\cdot\jump{\kappa^{\frac{1}{2}}\ph}\ds
		\\&\geq{C}\DGnorm{(\uh,\ph)}^{2}.
	\end{align*}
	The proof is complete.
\end{proof}
\begin{theorem}\label{thm:DG-error-estimate}
	Let $(\ud,\pd)\in[{H}^{\mv+1}(\od)]^{d}\times\left({H}^{\mq+1}(\od)\cap{L}^{2}_{0}(\od)\right)$ and $(\uh,\ph)\in\Uh\times\Ph$ be the solutions to problems~\eqref{eq:Darcy-Model} and~\eqref{eq:Darcy-discrete-problem}, respectively. Then the following error estimate holds:
	\begin{equation}\label{ineq:DG-error-estimate}
		\DGnorm{(\ud-\uh,\pd-\ph)}\leq{C}h^{s}(\|\ud\|_{H^{\mv+1}(\od)}+\|p\|_{H^{\mq+1}(\od)}),
	\end{equation}
	where $s=\min\{\mv,\mq\}$.
\end{theorem}
\begin{proof}
	First, we split the error as $\ud-\uh=(\ud-\Rv\ud)+(\Rv\ud-\uh)=\boldsymbol{\eta}_{\ud}+\boldsymbol{\theta}_{\ud}$, and $\pd-\ph=(\pd-\Rq\pd)+(\Rq\pd-\ph)=\eta_{\pd}+\theta_{\pd}$. From the error equation~\eqref{eq:error-equation} and~\eqref{ineq:coercive}, we obtain
	\begin{align*}
		\DGnorm{(\boldsymbol{\theta}_{\ud},\theta_{p})}^{2}&\leq{C}{\mathcal{A}}(\boldsymbol{\theta}_{\ud},\theta_{p};\boldsymbol{\theta}_{\ud},\theta_{p})=-{C}\mathcal{A}(\boldsymbol{\eta}_{\ud},\eta_{p};\boldsymbol{\theta}_{\ud},\theta_{p})\leq{C}\DGnorm{(\boldsymbol{\eta}_{\ud},\eta_{p})}\DGnorm{(\boldsymbol{\theta}_{\ud},\theta_{p})}.
	\end{align*}
	Using the properties of the reconstruction operator, we have
	\begin{align*}
		\DGnorm{(\boldsymbol{\eta}_{\ud},\eta_{p})}\leq{C}\big(h^{\mv}\|\ud\|_{H^{\mv+1}}+h^{\mq}\|\pd\|_{H^{\mq+1}}\big).
	\end{align*}
	Combining the above inequality and applying the triangle inequality yields the estimate~\eqref{ineq:DG-error-estimate}.
\end{proof}

Since the pressure error $\pd-\ph$ has zero mean, the discrete Poincar\'e inequality gives $\|\ud-\uh\|_{L^{2}(\od)}+\|\pd-\ph\|_{L^{2}(\od)}\leq{C}\DGnorm{(\ud-\uh,\pd-\ph)}$. Combining this inequality with the preceding energy norm estimate, we obtain the following $L^2$-error estimate.
\begin{theorem}\label{thm:L2-error-estimate}
	Let $(\ud,\pd)\in[{H}^{\mv+1}(\od)]^{d}\times\left({H}^{\mq+1}(\od)\cap{L}^{2}_{0}(\od)\right)$ and $(\uh,\ph)\in\Uh\times\Ph$ be the solutions to problems~\eqref{eq:Darcy-Model} and~\eqref{eq:Darcy-discrete-problem}, respectively. Then the following error estimate holds:
	\begin{equation}\label{ineq:L2-error-estimate}
		\|\ud-\uh\|_{L^{2}(\od)}+\|\pd-\ph\|_{L^{2}(\od)}\leq{C}h^{s}\big(\|\ud\|_{H^{\mv+1}(\od)}+\|\pd\|_{H^{\mq+1}(\od)}\big),
	\end{equation}
	where $s=\min\{\mv,\mq\}$.
\end{theorem}

	% vim:spell:tw=70:fo+=Mn:cc=70
\section{Preconditioning Method}\label{sec:Precondition}

Although the previous section establishes the discrete formulation and its error estimates, efficiently solving the linear system remains a major challenge. To this end, we introduce a preconditioning method in this section.

Let $\mathbf{u},\mathbf{p}$ denote the coefficient vectors of the velocity and pressure, respectively. Define
\begin{align*}
	\Mat:=\begin{bmatrix}M&B\\-B^{\trans}&S\end{bmatrix},\quad\bfx:=\begin{bmatrix}\mathbf{u}\\\mathbf{p}\end{bmatrix},\quad\mathbf{f}:=\begin{bmatrix}F_{v}\\F_{q}\end{bmatrix},
\end{align*}
where $M\in\mathbb{R}^{d\dimTh\times{d\dimTh}}$, $B\in\mathbb{R}^{d\dimTh\times{(\dimTh-1)}}$, $S\in\mathbb{R}^{(\dimTh-1)\times{(\dimTh-1)}}$ are the matrices associated with the corresponding bilinear forms, while $F_{v}$ and $F_{q}$ are the vectors associated with the corresponding linear functionals. Then the discrete formulation~\eqref{eq:Darcy-discrete-problem} yields the following linear system:
\begin{equation}\label{eq:linear-system}
	\Mat\bfx=\mathbf{f}.
\end{equation}
Since the elementwise gradients vanish on the piecewise constant space and the discrete DG energy norm separates into velocity and pressure contributions, we introduce the following bilinear form to construct the corresponding preconditioner:
\begin{align*}
	\mathcal{A}^{0}(\uh,\ph;\vh,\qh)=&\,\sum\limits_{K\in\Th}\int_{K}\xi\uh\cdot\vh\dx+\sum\limits_{e\in\ehi}\int_{e}\frac{\alpha}{h_{e}}\jump{\xi^{\frac{1}{2}}\uh\otimes\nd}:\jump{\xi^{\frac{1}{2}}\vh\otimes\nd}\ds\\&\,+\sum\limits_{e\in\ehb}\int_{e}\frac{\alpha}{h_{e}}\jump{\xi^{\frac{1}{2}}\uh\cdot\nd}\jump{\xi^{\frac{1}{2}}\vh\cdot\nd}\ds+\sum\limits_{e\in\ehi}\int_{e}\frac{\beta}{h_{e}}\jump{\kappa^{\frac{1}{2}}\ph}\cdot\jump{\kappa^{\frac{1}{2}}\qh}\ds,\\&\,\forall\,\uh,\vh\in\Uconst,\;p_{h},q_{h}\in\Pmean.
\end{align*}
Let $\MatPre$ be the matrix associated with the bilinear form $\mathcal{A}^{0}$. To analyze the efficiency of the preconditioned GMRES solver, we first present the following key norm-equivalence lemma.
\begin{lemma}\label{lem:norm-equivalence-DG}
	Assume that the parameters $\alpha$ and $\beta$ are positive. Then there exists a constant $C>0$, depending on $\alpha$ and $\beta$ but independent of $h$, such that for any $\uh\in\Uconst$ and $\ph\in\Pmean$,
	\begin{equation}
		\DGnorm{(\uh,\ph)}\leq{C}\DGnorm{(\Rv\uh,\Rq\ph)}\leq{C}\overline{\Lambda}_{m}\DGnorm{(\uh,\ph)},
	\end{equation}
	where $\overline{\Lambda}_{m}=\max\{\overline{\Lambda}_{\mq},\overline{\Lambda}_{\mv}\}$, with $\overline{\Lambda}_{\mq}$ and $\overline{\Lambda}_{\mv}$ defined in Lemma~\ref{lem:Recons-error-estimate-scalar} and Lemma~\ref{lem:Recons-error-estimate-vector}, respectively.
\end{lemma}
\begin{proof}
	First, for any interior face $e\in\ehi$ and element $K\in\Th$, we have
	\begin{align*}
		\int_{e}\frac{\beta}{h_{e}}\jump{\kappa^{\frac{1}{2}}\Rq\ph}\cdot\jump{\kappa^{\frac{1}{2}}\Rq\ph}\ds&=\int_{e}\frac{\beta}{h_{e}}\jump{\kappa^{\frac{1}{2}}\Rqsca\ph}\cdot\jump{\kappa^{\frac{1}{2}}\Rqsca\ph}\ds,\\
		\int_{K}\kappa(\nabla\Rq\ph)\cdot(\nabla\Rq\ph)\dx&=\int_{K}\kappa(\nabla\Rqsca\ph)\cdot(\nabla\Rqsca\ph)\dx.
	\end{align*}
	Thus, the equivalence between
	\begin{align*}
		\DGnorm{\ph}\leq{C}\DGnorm{\Rq\ph}\leq{C}\overline{\Lambda}_{m}\DGnorm{\ph}
	\end{align*}
	and
	\begin{align}\label{ineq:norm-eq-p}
		\DGnorm{\ph}\leq{C}\DGnorm{\Rqsca\ph}\leq{C}\overline{\Lambda}_{m}\DGnorm{\ph}
	\end{align}
	follows. The estimate in~\eqref{ineq:norm-eq-p} follows directly from~\cite[Lemma~7]{LiRuo2024Preconditioned}, while the following velocity estimate is obtained by applying the same argument componentwise:
	\begin{align*}
		&\quad\sum\limits_{e\in\ehi}h_{e}^{-1}\Ltnorm{\jump{\uh\otimes\nd}}{e}^{2}+\sum\limits_{e\in\ehb}h_{e}^{-1}\Ltnorm{\jump{\uh\cdot\nd}}{e}^{2}\\&\leq{C}\bigg(\sum\limits_{K\in\Th}\Ltnorm{\nabla\Rv\uh}{K}^{2}+\sum\limits_{e\in\ehi}h_{e}^{-1}\Ltnorm{\jump{\Rv\uh\otimes\nd}}{e}^{2}+\sum\limits_{e\in\ehb}h_{e}^{-1}\Ltnorm{\jump{\Rv\uh\cdot\nd}}{e}^{2}\bigg)\\&\leq{C}\overline{\Lambda}_{m}^{2}\bigg(\sum\limits_{e\in\ehi}h_{e}^{-1}\Ltnorm{\jump{\uh\otimes\nd}}{e}^{2}+\sum\limits_{e\in\ehb}h_{e}^{-1}\Ltnorm{\jump{\uh\cdot\nd}}{e}^{2}\bigg).
	\end{align*}
	The details are omitted. Hence, it suffices to prove the equivalence for
	\begin{align*}
		\sum\limits_{K\in\Th}\Ltnorm{\xi^{\frac{1}{2}}\uh}{K}^{2}\leq{C}\sum\limits_{K\in\Th}\Ltnorm{\xi^{\frac{1}{2}}\Rv\uh}{K}^{2}\leq{C}\overline{\Lambda}_{m}^{2}\sum\limits_{K\in\Th}\Ltnorm{\xi^{\frac{1}{2}}\uh}{K}^{2}.
	\end{align*}
	
	We first prove the above inequality starting with the lower bound. For any $K\in\Th$, 
	\begin{align*}
		(\Rv\uhi)(\x_{K})=(\Rv_{K}\uhi)(\x_{K})=\uhi(\x_{K}),
	\end{align*}
	where $\uhi$ denotes the $i$-th component of $\uh$. This identity implies that
	\begin{align*}
		\|\xi^{\frac{1}{2}}\uhi\|_{L^{2}(K)}^{2}&\leq{C}\xi_{\max}h_{K}^{d}|\uhi(\x_{K})|^{2}={C}h_{K}^{d}|\Rv\uhi(\x_{K})|^{2}\leq{C}h_{K}^{d}\|\Rv\uhi\|_{L^{\infty}(K)}^{2}\\&\leq{C}\|\xi^{\frac{1}{2}}\Rv_{K}\uhi\|_{L^{2}(K)}^{2}.
	\end{align*}
	The lower bound estimate follows by summing over all elements and all components. Regarding the upper bound, for any $K\in\Th$, we derive that
	\begin{align*}
		\|\xi^{\frac{1}{2}}\Rv\uhi\|_{L^{2}(K)}^{2}\leq{C}\xi_{\max}h_{K}^{d}\|\Rv\uhi\|_{L^{\infty}(K)}^{2}\leq{C}\overline{\Lambda}_{m}^{2}\max\limits_{\x\in{I}(K)}h_{K}^{d}|\uhi(\x)|^{2}.
	\end{align*}
	Let $K'\in{S(K)}$ be such that $\max_{\x\in{I}(K)}|\uhi(\x)|=|\uhi(\x_{K'})|$. Similarly, summing over all elements and using the finite overlap property of the element patches, we have
	\begin{align*}
		\sum\limits_{K\in\Th}\|\xi^{\frac{1}{2}}\Rv\uh\|_{L^{2}(K)}^{2}&\leq{C}\overline{\Lambda}_{m}^{2}\sum\limits_{i=1}^{d}\sum\limits_{K\in\Th}h_{K}^{d}\|\uhi\|_{L^{\infty}(S(K))}^{2}\leq{C}\overline{\Lambda}_{m}^{2}\sum\limits_{i=1}^{d}\sum\limits_{K\in\Th}\|\xi^{\frac{1}{2}}\uhi\|_{L^{2}(K')}^{2}
		\\&\leq{C}\overline{\Lambda}_{m}^{2}\sum\limits_{K'\in\Th}\|\xi^{\frac{1}{2}}\uh\|_{L^{2}(K')}^{2}.
	\end{align*}
	The proof is complete.
\end{proof}

Let $\mathbf{u}_{h},\mathbf{v}_{h}\in\mathbb{R}^{d\dimTh}$ and $\mathbf{p}_{h},\mathbf{q}_{h}\in\mathbb{R}^{\dimTh-1}$ denote the coefficient vectors corresponding to $\uh,\vh\in\Uconst$ and $\ph,\qh\in\Pmean$, respectively. We set $\bfx:=(\mathbf{u}_{h},\mathbf{p}_{h})^{\trans}$ and $\bfy:=(\mathbf{v}_{h},\mathbf{q}_{h})^{\trans}$. Then, by Lemma~\ref{lem:continuity} and Lemma~\ref{lem:coercive}, we have
\begin{align*}
	|\bfy^{\trans}\Mat\bfx|&=\left|\mathcal{A}(\Rv\uh,\Rq\ph;\Rv\vh,\Rq\qh)\right|\leq{C}\DGnorm{(\Rv\uh,\Rq\ph)}\DGnorm{(\Rv\vh,\Rq\qh)}
	\\&\leq{C}\overline{\Lambda}_{m}^{2}\DGnorm{(\uh,\ph)}\DGnorm{(\vh,\qh)}={C}\overline{\Lambda}_{m}^{2}\sqrt{\mathcal{A}^{0}(\uh,\ph;\uh,\ph)}\sqrt{\mathcal{A}^{0}(\vh,\qh;\vh,\qh)}
	\\&\,={C}\overline{\Lambda}_{m}^{2}\sqrt{\bfx^{\trans}\MatPre\bfx}\sqrt{\bfy^{\trans}\MatPre\bfy},\\
	\bfx^{\trans}\Mat\bfx&=\mathcal{A}(\Rv\uh,\Rq\ph;\Rv\uh,\Rq\ph)\geq{C}\DGnorm{(\Rv\uh,\Rq\ph)}^{2}
	\\&\geq{c}\DGnorm{(\uh,\ph)}^{2}={c}\mathcal{A}^{0}(\uh,\ph;\uh,\ph)
	\\&\,={c}\bfx^{\trans}\MatPre\bfx.
\end{align*}

Since $\mathcal{A}^{0}$ is symmetric and positive definite on $\Uconst\times\Pmean$, the same holds for the preconditioning matrix $\MatPre$. Consequently, $(\bfx,\bfy)_{\MatPre}:=\bfx^{\trans}\MatPre\bfy$ defines an inner product with the induced norm $\|\bfx\|_{\MatPre}:=
(\bfx^{\trans}\MatPre\bfx)^{1/2}$, and the preconditioned matrix $\MatPre^{-1}\Mat$ is well defined on $\mathbb{R}^{d\dimTh+(\dimTh-1)}$. To establish the effectiveness of the preconditioner for GMRES, we follow the analysis in~\cite{XuJinchao1992GMRES}. It suffices to prove that there exist positive constants $c$ and $C$ such that
\begin{equation}\label{ineq:sv-sup}
	\|\MatPre^{-1}\Mat\bfx\|_{\MatPre}^{2}\leq{C}\overline{\Lambda}_{m}^{4}\|\bfx\|_{\MatPre}^{2},
\end{equation}
and
\begin{equation}\label{ineq:ev-inf}
	(\MatPre^{-1}\Mat\bfx,\bfx)_{\MatPre}\geq{c}\|\bfx\|_{\MatPre}^{2}.
\end{equation}
We first prove the lower bound. For any $\bfx\in\mathbb{R}^{d\dimTh+(\dimTh-1)}$,
\begin{align*}
	(\MatPre^{-1}\Mat\bfx,\bfx)_{\MatPre}=\bfx^{\trans}\Mat\bfx\geq{c}\bfx^{\trans}\MatPre\bfx={c}\|\bfx\|_{\MatPre}^{2}.
\end{align*}
Next, we prove the upper bound. For any $\bfx\in\mathbb{R}^{d\dimTh+(\dimTh-1)}$, we obtain
\begin{align*}
	\|\MatPre^{-1}\Mat\bfx\|_{\MatPre}=\sup\limits_{\substack{\bfy\in\mathbb{R}^{d\dimTh+(\dimTh-1)}\\\bfy\neq\mathbf{0}}}\frac{|\bfy^{\trans}\Mat\bfx|}{\|\bfy\|_{\MatPre}}\leq{C}\overline{\Lambda}_{m}^{2}\|\bfx\|_{\MatPre}.
\end{align*}

Based on the above two estimates, we apply~\cite[Theorem~2]{XuJinchao1992GMRES} to obtain the following theorem.
\begin{theorem}\label{thm:GMRES-rates}
	There exist positive constants $c$ and $C$, independent of $h$, such that the preconditioned GMRES method converges and satisfies
	\begin{equation}\label{eq:GMRES-rates}
		\|\mathbf{r}_{k}\|_{\MatPre}\leq\bigg(1-\bigg(\frac{c}{C\overline{\Lambda}_{m}^{2}}\bigg)^{2}\bigg)^{\frac{k}{2}}\|\mathbf{r}_{0}\|_{\MatPre},
	\end{equation}
	where $\mathbf{r}_{k}=\MatPre^{-1}(\mathbf{f}-\Mat\bfx_{k})$, and $\bfx_{k}$ denotes the $k$-th GMRES iterate.
\end{theorem}
Consequently, the $\MatPre$-preconditioned GMRES method exhibits mesh-independent convergence. However, since the corresponding convergence factor depends on $\overline{\Lambda}_{m}$, the convergence rate is influenced by the polynomial degrees $\mv$ and $\mq$, as well as the patch cardinality $\#S(K)$.

In the implementation, the Lagrange-multiplier formulation is used to impose the zero-mean pressure constraint~\cite{Bochev2005Neumann}. Accordingly, the preconditioner is applied by solving the velocity block system and the augmented pressure-constraint system. For large-scale discrete systems, direct block solvers can be computationally expensive and require substantial memory. Therefore, multigrid methods can be employed to approximate the actions of the block inverses. In this work, we employ GMG and BoomerAMG provided by the HYPRE package~\cite{Falgout2006design}.

	% vim:spell:tw=70:fo+=Mn:cc=70
\section{Numerical Experiments}\label{sec:Experiments}
In this section, several numerical examples are presented to assess the convergence behavior of the numerical scheme and the effectiveness of the proposed preconditioner for GMRES.

\subsection{Smooth Problem}
First, we consider two- and three-dimensional problems with smooth exact solutions, discretized using triangular and tetrahedral meshes, respectively. With $\kappa=1$, the corresponding source term $f$ and boundary data $g$ are determined by substituting the prescribed exact solutions into the governing equations.
\begin{example}\label{ex:2D-smooth}
	In this example, the computational domain is $\od=(0,1)\times(0,1)$. The exact solutions for the velocity and pressure are given by
	\begin{align*}
		\ud(x,y)&=\begin{pmatrix}-\pi\sin(2\pi{x})\sin(\pi{y})^{2}\\-\pi\sin(2\pi{y})\sin(\pi{x})^{2}\end{pmatrix},\\\pd(x,y)&=\sin(\pi{x})^{2}\sin(\pi{y})^{2}-\frac{1}{4}.
	\end{align*}
\end{example}
\begin{example}\label{ex:3D-smooth}
	In this example, the computational domain is $\od=(0,1)^{3}$. The exact solutions for the velocity and pressure are given by
	\begin{align*}
		\ud(x,y,z)&=\begin{pmatrix}-\pi\cos(\pi{x})\sin(\pi{y})\sin(\pi{z})\\-\pi\sin(\pi{x})\cos(\pi{y})\sin(\pi{z})\\-\pi\sin(\pi{x})\sin(\pi{y})\cos(\pi{z})\end{pmatrix},
		\\\pd(x,y,z)&=\sin(\pi{x})\sin(\pi{y})\sin(\pi{z})-\frac{8}{\pi^{3}}.
	\end{align*}
\end{example}

To simplify the notation, we set $e_{\ud}=\Ltnorm{\ud-\uh}{\od}$, $e_{\pd}=\Ltnorm{\pd-\ph}{\od}$, and $e_{\mathrm{DG}}=\DGnorm{(\ud-\uh,\pd-\ph)}$.
\begin{table}[!t]
	\caption{The convergence orders for Example~\ref{ex:2D-smooth}.}
	\label{tab:2D-smooth-order}
	\centering
	\renewcommand{\arraystretch}{1.2}
	\begin{tabular}{cc*{6}{>{\centering\arraybackslash}p{4.5em}}c}
		\hline
		$m$&${1}/{h}$&$10$&$20$&$40$&$80$&$160$&$320$&Order\\\hline
		\multirow{3}{*}{$1$}
		&${e}_{\ud}$&$1.58\mathrm{e}{-01}$&$5.77\mathrm{e}{-02}$&$1.33\mathrm{e}{-02}$&$3.36\mathrm{e}{-03}$&$8.99\mathrm{e}{-04}$&$2.36\mathrm{e}{-04}$&1.93\\
		&$e_{\pd}$&$3.06\mathrm{e}{-02}$&$1.29\mathrm{e}{-02}$&$2.91\mathrm{e}{-03}$&$7.53\mathrm{e}{-04}$&$1.69\mathrm{e}{-04}$&$4.44\mathrm{e}{-05}$&1.93\\
		&$e_{\mathrm{DG}}$&$3.80\mathrm{e}{+00}$&$1.99\mathrm{e}{+00}$&$1.00\mathrm{e}{+00}$&$5.07\mathrm{e}{-01}$&$2.55\mathrm{e}{-01}$&$1.28\mathrm{e}{-01}$&0.99\\\hline
		\multirow{3}{*}{$2$}
		&${e}_{\ud}$&$6.77\mathrm{e}{-02}$&$2.01\mathrm{e}{-02}$&$5.39\mathrm{e}{-03}$&$1.39\mathrm{e}{-03}$&$3.51\mathrm{e}{-04}$&$8.82\mathrm{e}{-05}$&1.99\\
		&$e_{\pd}$&$4.66\mathrm{e}{-03}$&$9.85\mathrm{e}{-04}$&$2.52\mathrm{e}{-04}$&$6.28\mathrm{e}{-05}$&$1.56\mathrm{e}{-05}$&$3.92\mathrm{e}{-06}$&2.00\\
		&$e_{\mathrm{DG}}$&$8.25\mathrm{e}{-01}$&$2.10\mathrm{e}{-01}$&$5.23\mathrm{e}{-02}$&$1.31\mathrm{e}{-02}$&$3.27\mathrm{e}{-03}$&$8.18\mathrm{e}{-04}$&2.00\\\hline
		\multirow{3}{*}{$3$}
		&${e}_{\ud}$&$9.68\mathrm{e}{-03}$&$1.08\mathrm{e}{-03}$&$7.65\mathrm{e}{-05}$&$5.00\mathrm{e}{-06}$&$3.26\mathrm{e}{-07}$&$2.08\mathrm{e}{-08}$&3.97\\
		&$e_{\pd}$&$1.25\mathrm{e}{-03}$&$2.56\mathrm{e}{-04}$&$1.90\mathrm{e}{-05}$&$1.22\mathrm{e}{-06}$&$7.89\mathrm{e}{-08}$&$5.05\mathrm{e}{-09}$&3.96\\
		&$e_{\mathrm{DG}}$&$2.68\mathrm{e}{-01}$&$3.25\mathrm{e}{-02}$&$4.09\mathrm{e}{-03}$&$5.16\mathrm{e}{-04}$&$6.53\mathrm{e}{-05}$&$8.23\mathrm{e}{-06}$&2.99\\\hline
		\multirow{3}{*}{$4$}
		&${e}_{\ud}$&$4.09\mathrm{e}{-03}$&$2.70\mathrm{e}{-04}$&$1.91\mathrm{e}{-05}$&$1.29\mathrm{e}{-06}$&$8.41\mathrm{e}{-08}$&$5.37\mathrm{e}{-09}$&3.97\\
		&$e_{\pd}$&$2.77\mathrm{e}{-04}$&$1.26\mathrm{e}{-05}$&$7.28\mathrm{e}{-07}$&$4.79\mathrm{e}{-08}$&$3.20\mathrm{e}{-09}$&$2.12\mathrm{e}{-10}$&3.92\\
		&$e_{\mathrm{DG}}$&$7.39\mathrm{e}{-02}$&$3.85\mathrm{e}{-03}$&$2.25\mathrm{e}{-04}$&$1.35\mathrm{e}{-05}$&$8.19\mathrm{e}{-07}$&$5.04\mathrm{e}{-08}$&4.02\\\hline
	\end{tabular}
\end{table}
\begin{table}[!t]
	\caption{The convergence orders for Example~\ref{ex:3D-smooth}.}
	\label{tab:3D-smooth-order}
	\centering
	\renewcommand{\arraystretch}{1.2}
	\begin{tabular}{cc*{4}{>{\centering\arraybackslash}p{4.5em}}c}
		\hline
		$m$&${1}/{h}$&$4$&$8$&$16$&$32$&Order\\\hline
		\multirow{3}{*}{$1$}
		&${e}_{\ud}$&$3.06\mathrm{e}{-01}$&$7.65\mathrm{e}{-02}$&$1.81\mathrm{e}{-02}$&$4.31\mathrm{e}{-03}$&2.07\\
		&$e_{\pd}$&$3.86\mathrm{e}{-02}$&$1.11\mathrm{e}{-02}$&$3.21\mathrm{e}{-03}$&$8.62\mathrm{e}{-04}$&1.90\\
		&$e_{\mathrm{DG}}$&$5.10\mathrm{e}{+00}$&$2.47\mathrm{e}{+00}$&$1.21\mathrm{e}{+00}$&$6.00\mathrm{e}{-01}$&1.01\\\hline
		\multirow{3}{*}{$2$}
		&${e}_{\ud}$&$1.21\mathrm{e}{-01}$&$2.75\mathrm{e}{-02}$&$6.39\mathrm{e}{-03}$&$1.51\mathrm{e}{-03}$&2.08\\
		&$e_{\pd}$&$1.28\mathrm{e}{-02}$&$2.03\mathrm{e}{-03}$&$4.69\mathrm{e}{-04}$&$1.16\mathrm{e}{-04}$&2.01\\
		&$e_{\mathrm{DG}}$&$1.65\mathrm{e}{+00}$&$3.89\mathrm{e}{-01}$&$9.68\mathrm{e}{-02}$&$2.43\mathrm{e}{-02}$&2.00\\\hline
		\multirow{3}{*}{$3$}
		&${e}_{\ud}$&$2.70\mathrm{e}{-02}$&$1.98\mathrm{e}{-03}$&$1.43\mathrm{e}{-04}$&$9.75\mathrm{e}{-06}$&3.88\\
		&$e_{\pd}$&$4.87\mathrm{e}{-03}$&$3.32\mathrm{e}{-04}$&$2.23\mathrm{e}{-05}$&$1.50\mathrm{e}{-06}$&3.89\\
		&$e_{\mathrm{DG}}$&$5.60\mathrm{e}{-01}$&$6.34\mathrm{e}{-02}$&$7.36\mathrm{e}{-03}$&$8.82\mathrm{e}{-04}$&3.06
		\\\hline
	\end{tabular}
\end{table}

In Table~\ref{tab:2D-smooth-order} and Table~\ref{tab:3D-smooth-order}, we report the velocity and pressure errors in the $L^2$ norm and the DG energy norm for the equal-order case $m_{\boldsymbol{v}}=m_{q}=m$. As observed, the errors $e_{\ud}$, $e_{\pd}$, and $e_{\mathrm{DG}}$ converge at the rates $\mathcal{O}(h^{m+1})$, $\mathcal{O}(h^{m+1})$, and $\mathcal{O}(h^{m})$, respectively, for odd $m$, whereas for even $m$, all three errors converge at the rate $\mathcal{O}(h^{m})$. The numerical results confirm the validity of Theorem~\ref{thm:DG-error-estimate} and Theorem~\ref{thm:L2-error-estimate}. Although the theoretical analysis yields only an $\mathcal{O}(h^{m})$ estimate, the numerical results indicate that the $L^{2}$ errors may attain the optimal
convergence rate $\mathcal{O}(h^{m+1})$ for odd polynomial degrees. 

To further investigate the approximation performance of the proposed method, we compare the numerical errors of the RDA method with those of the DG method for different numbers of degrees of freedom.

\begin{figure}[!t]
	\centering
	\begin{subfigure}[t]{0.48\textwidth}
		\centering
		\includegraphics[width=0.75\textwidth]{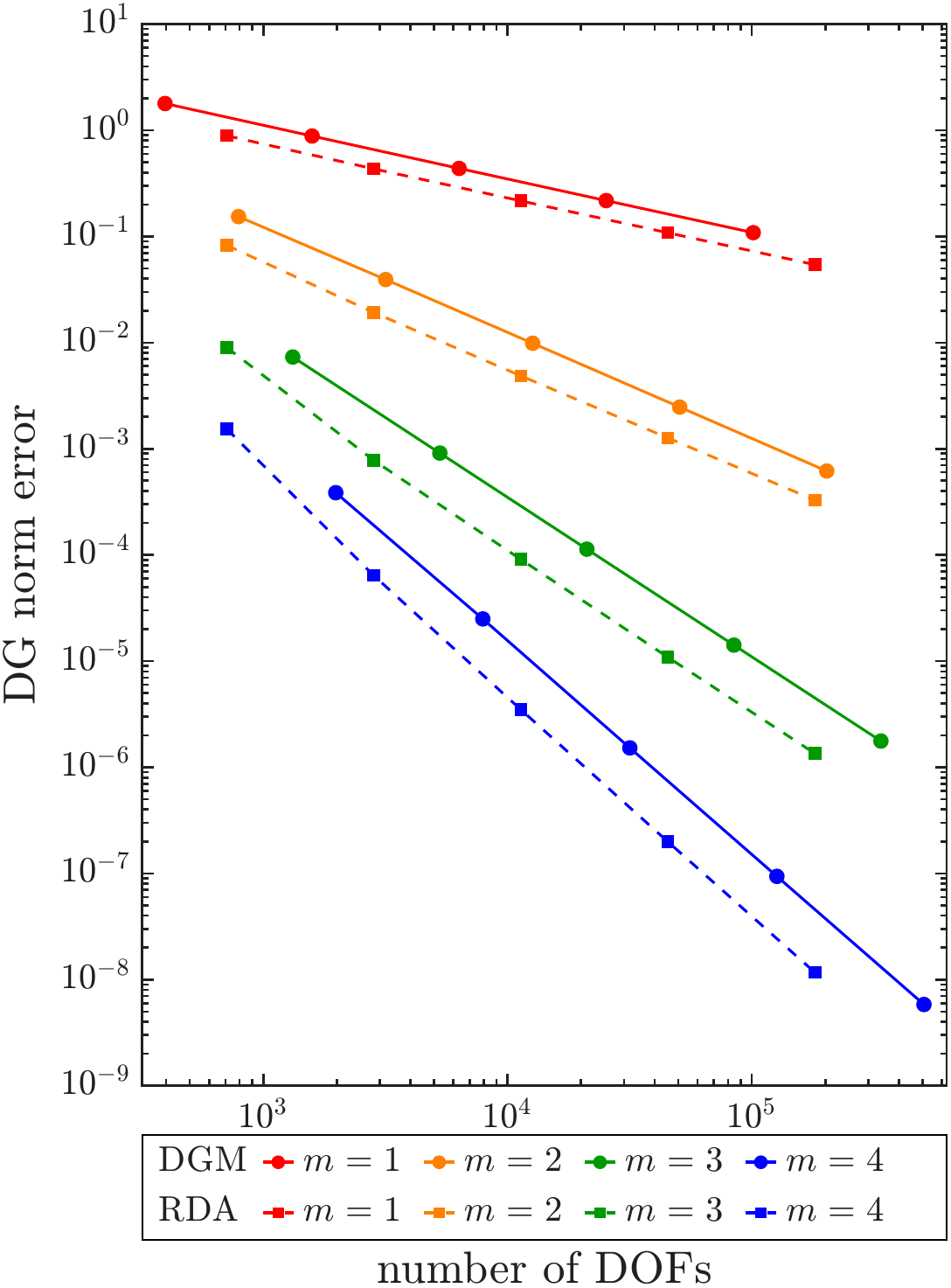}
	\end{subfigure}
	\begin{subfigure}[t]{0.48\textwidth}
		\centering
		\includegraphics[width=0.75\textwidth]{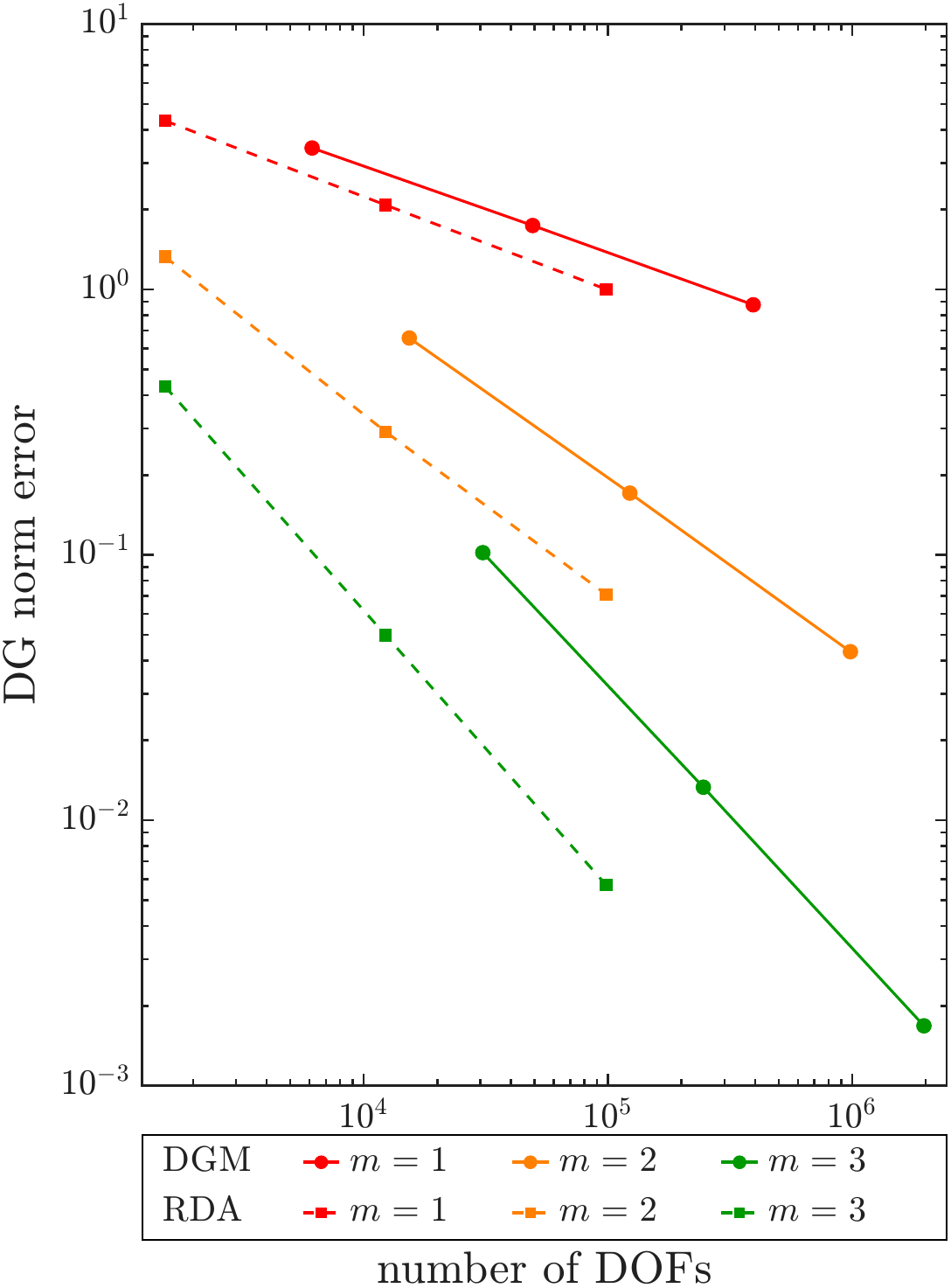}
	\end{subfigure}
	\caption{Comparison of DG norm errors for the RDA and standard DG methods for Example~\ref{ex:2D-smooth} (left) and Example~\ref{ex:3D-smooth} (right).}
	\label{fig:RDA-DG}
\end{figure}

From Figure~\ref{fig:RDA-DG}, both the RDA and standard DG methods exhibit $m$-th order convergence in the DG norm. For the same number of degrees of freedom, RDA generally yields smaller errors, indicating the effectiveness of the reconstruction in improving the approximation accuracy.

We next investigate the computational efficiency of solving the discrete systems and assess the effectiveness of the proposed preconditioned solver. The numerical experiments compare two approaches: unpreconditioned GMRES, and $\MatPre$-preconditioned GMRES. For the proposed preconditioner, we further consider direct and approximate block solvers based on GMG and HYPRE BoomerAMG. The different approaches are compared in terms of iteration counts and CPU time. We first present the spectral properties of the preconditioned discrete system. 

Table~\ref{tab:2D-smooth-spectral} and Table~\ref{tab:3D-smooth-spectral} show that, for a fixed approximation degree $m$, the maximum singular value $\sigma_{\max}$ of the matrix $\MatPre^{-1/2}\Mat\MatPre^{-1/2}$ and the minimum eigenvalue $\lambda_{\min}$ of its symmetric part remain nearly unchanged under mesh refinement. These results are consistent with the mesh-independent field-of-values bounds established for the preconditioned system.

\begin{table}[!t]
	\caption{Spectral properties of the preconditioned discrete system for Example~\ref{ex:2D-smooth}.}
	\label{tab:2D-smooth-spectral}
	\renewcommand{\arraystretch}{1.2}
	\centering
	\begin{tabular}{ccccccccc}
		\hline
		$m$&\multicolumn{2}{c}{$1$}&\multicolumn{2}{c}{$2$}&\multicolumn{2}{c}{$3$}&\multicolumn{2}{c}{$4$}\\
		${1}/{h}$&$\sigma_{\max}$&$\lambda_{\min}$&$\sigma_{\max}$&$\lambda_{\min}$&$\sigma_{\max}$&$\lambda_{\min}$&$\sigma_{\max}$&$\lambda_{\min}$\\\hline
		$10$ &3.21&0.45&3.69&0.38&6.44&0.37&13.21&0.37\\
		$20$ &3.34&0.41&3.78&0.36&7.01&0.34&17.35&0.34\\
		$40$ &3.25&0.39&4.03&0.35&6.92&0.33&16.82&0.34\\
		$80$ &3.26&0.39&4.34&0.34&6.94&0.32&16.77&0.33\\
		$160$&3.27&0.39&4.41&0.34&6.96&0.32&16.77&0.33
		\\\hline
	\end{tabular}
\end{table}

\begin{table}[!t]
	\caption{Spectral properties of the preconditioned discrete system for Example~\ref{ex:3D-smooth}.}
	\label{tab:3D-smooth-spectral}
	\renewcommand{\arraystretch}{1.2}
	\centering
	\begin{tabular}{ccccccc}
		\hline
		$m$&\multicolumn{2}{c}{$1$}&\multicolumn{2}{c}{$2$}&\multicolumn{2}{c}{$3$}\\
		${1}/{h}$&$\sigma_{\max}$&$\lambda_{\min}$&$\sigma_{\max}$&$\lambda_{\min}$&$\sigma_{\max}$&$\lambda_{\min}$\\\hline
		$4$ &3.49&0.33&4.81&0.28&12.51&0.27\\
		$8$ &3.58&0.34&4.88&0.28&12.79&0.27\\
		$16$&3.81&0.33&4.86&0.28&12.75&0.27
		\\\hline
	\end{tabular}
\end{table}
\begin{table}[!t]
	\caption{Iteration counts of different iterative solvers in Example~\ref{ex:2D-smooth}.}
	\label{tab:2D-smooth-steps}
	\renewcommand{\arraystretch}{1.2}
	\centering
	\begin{threeparttable}
		\begin{tabular}{ccc*{6}{>{\centering\arraybackslash}p{3em}}}
			\hline
			\multirow{2}{*}{$m$}&\multirow{2}{*}{Method}
			& \multirow{2}{*}{Block solver}&\multicolumn{6}{c}{$1/h$}\\\cline{4-9}
			& & & 10 & 20 & 40 & 80 & 160 & 320 \\\hline
			\multirow{4}{*}{$1$}&
			GMRES&--&183 &306 &620 &1312&2975&NC\\
			&\multirow{3}{*}{$\MatPre$-GMRES}&Direct&24&23&22&21&20&20\\
			&&GMG&24&23&22&21&21&21\\
			&&BoomerAMG&26&26&24&23&23&23\\\hline
			\multirow{4}{*}{$2$}&
			GMRES&--&176 &335 &643 &1605&3791&NC\\
			&\multirow{3}{*}{$\MatPre$-GMRES}&Direct&26&25&24&24&22&21\\
			&&GMG&26&25&25&24&24&23\\
			&&BoomerAMG&28&28&27&27&26&25\\\hline
			\multirow{4}{*}{$3$}&
			GMRES&--&229 &443 &853 &2448&NC&NC\\
			&\multirow{3}{*}{$\MatPre$-GMRES}&Direct&42&45&43&40&39&37\\
			&&GMG&42&45&43&41&39&38\\
			&&BoomerAMG&45&48&46&44&42&40\\\hline
			\multirow{4}{*}{$4$}&
			GMRES&--&251 &492 &950 &3774&NC&NC\\
			&\multirow{3}{*}{$\MatPre$-GMRES}&Direct&67&75&77&75&72&70\\
			&&GMG&67&75&78&76&74&73\\
			&&BoomerAMG&71&81&82&80&78&75
			\\\hline
		\end{tabular}
		\begin{tablenotes}
			\footnotesize
			\item \makebox[4em][l]{NC} indicates that the iterative method fails to converge within 5000 iterations.
		\end{tablenotes}
	\end{threeparttable}
\end{table}
\begin{table}[!t]
	\caption{Iteration counts of different iterative solvers in Example~\ref{ex:3D-smooth}.}
	\label{tab:3D-smooth-steps}
	\renewcommand{\arraystretch}{1.2}
	\centering
	\begin{tabular}{ccc*{4}{>{\centering\arraybackslash}p{3em}}}
		\hline
		\multirow{2}{*}{$m$}&\multirow{2}{*}{Method}
		& \multirow{2}{*}{Block solver}&\multicolumn{4}{c}{$1/h$}\\\cline{4-7}
		& & & 4 & 8 & 16 & 32 \\\hline
		\multirow{4}{*}{$1$}&
		GMRES&--&115&203&362&659\\
		&\multirow{3}{*}{$\MatPre$-GMRES}&Direct&27&29&31&29\\
		&&GMG&27&31&31&31\\
		&&BoomerAMG&29&31&33&34\\\hline
		\multirow{4}{*}{$2$}&
		GMRES&--&128&213&370&604\\
		&\multirow{3}{*}{$\MatPre$-GMRES}&Direct&37&37&35&33\\
		&&GMG&37&38&37&36\\
		&&BoomerAMG&39&40&38&37\\\hline
		\multirow{4}{*}{$3$}&
		GMRES&--&165 &307 &579 &1037\\
		&\multirow{3}{*}{$\MatPre$-GMRES}&Direct&55&58&57&55\\
		&&GMG&55&59&59&58\\
		&&BoomerAMG&58&62&61&61
		\\\hline
	\end{tabular}
\end{table}

As shown in Table~\ref{tab:2D-smooth-steps} and Table~\ref{tab:3D-smooth-steps}, the iteration counts of unpreconditioned GMRES exhibit a strong dependence on the mesh size $h$. As the mesh is refined, the iteration count increases rapidly. In contrast, with the $\MatPre$-preconditioned GMRES method, all three block solvers achieve convergence within a relatively small number of iterations, and the iteration counts remain essentially stable under mesh refinement. This further confirms the efficient convergence and mesh robustness of the proposed preconditioner. On the other hand, the iteration counts still increase with the approximation degree $m$, consistent with Lemma~\ref{lem:norm-equivalence-DG}, where the estimates depend on the reconstruction constant $\overline{\Lambda}_m$. Therefore, the iteration counts of the $\MatPre$-preconditioned GMRES method are essentially independent of the mesh size $h$, while still depending on the approximation degree $m$.

\begin{table}[!t]
	\caption{CPU times (s) for different solution methods in Example~\ref{ex:2D-smooth}.}
	\label{tab:2D-smooth-times}
	\renewcommand{\arraystretch}{1.15}
	\centering
	\begin{threeparttable}
		\setlength{\tabcolsep}{0pt}
		\begin{tabular*}{\linewidth}{@{\extracolsep{\fill}}cc*{6}{c}@{}}
			\hline
			\multirow{3}{*}{$m$}&\multirow{3}{*}{${1}/{h}$}&\multirow{3}{*}{GMRES}&\multicolumn{5}{c}{$\MatPre$-GMRES}\\\cline{4-8}
			&&&\multirow{2}{*}{Assemble}&\multicolumn{2}{c}{Direct}&\multirow{2}{*}{GMG}&\multirow{2}{*}{BoomerAMG}\\\cline{5-6}
			&&&&Factor&Solve&&\\\hline
			\multirow{5}{*}{$1$}&20&0.04&0.01&--&--&--&--\\
			&40&0.87&0.02&0.01&0.01&0.01&0.01\\
			&80&10.45&0.11&0.03&0.06&0.05&0.05\\
			&160&109.85&0.52&0.19&0.33&0.24&0.24\\
			&320&NC&2.32&1.45&1.69&1.08&1.04\\\hline
			\multirow{5}{*}{$2$}&20&0.07&0.01&--&0.01&--&--\\
			&40&1.09&0.03&0.01&0.02&0.02&0.02\\
			&80&14.15&0.15&0.03&0.09&0.08&0.08\\
			&160&145.85&0.62&0.19&0.43&0.35&0.35\\
			&320&NC&2.58&1.44&1.98&1.50&1.38\\\hline
			\multirow{5}{*}{$3$}&20&0.12&0.01&--&0.01&0.01&0.01\\
			&40&2.03&0.03&0.01&0.04&0.04&0.05\\
			&80&24.46&0.12&0.03&0.19&0.18&0.19\\
			&160&NC&0.55&0.19&0.92&0.73&0.77\\
			&320&NC&2.14&1.44&4.13&3.21&3.02\\\hline
			\multirow{5}{*}{$4$}&20&0.18&0.01&--&0.02&0.02&0.03\\
			&40&2.70&0.03&0.01&0.11&0.11&0.11\\
			&80&43.88&0.15&0.03&0.49&0.45&0.45\\
			&160&NC&0.66&0.19&2.20&1.84&1.93\\
			&320&NC&2.69&1.45&9.52&7.85&7.36\\\hline
		\end{tabular*}
		\begin{tablenotes}
			\footnotesize
			\item \makebox[4em][l]{--} indicates that the CPU time is below the displayed precision.
			\item \makebox[4em][l]{NC} indicates that the iterative method fails to converge within 5000 iterations.
		\end{tablenotes}
	\end{threeparttable}
\end{table}

Table~\ref{tab:2D-smooth-times}~and~\ref{tab:3D-smooth-times} show that the solution time of unpreconditioned GMRES increases rapidly with mesh refinement. Consequently, solving the linear system without preconditioning becomes increasingly expensive on finer meshes. In contrast, the $\MatPre$-preconditioned GMRES method with direct block solvers significantly reduces the GMRES solve time compared with unpreconditioned GMRES on the finest three-dimensional meshes. However, the computational cost of direct block factorization remains substantial and can make the total solution time comparable to or even greater than that of unpreconditioned GMRES on fine three-dimensional meshes. To reduce this computational cost, GMG and BoomerAMG are employed as approximate block solvers. For GMG, the reported setup time includes the construction of the multigrid hierarchy and the assembly of the operator matrices on all levels. Compared with direct block solvers, these approximate solvers further reduce the GMRES solve time, with the improvement becoming more evident on the finest tested meshes. These results demonstrate that approximate block solvers can effectively improve the computational efficiency of the proposed preconditioned solver.

\begin{table}[!t]
	\caption{CPU times (s) for different solution methods in Example~\ref{ex:3D-smooth}.}
	\label{tab:3D-smooth-times}
	\renewcommand{\arraystretch}{1.2}
	\centering
	\begin{threeparttable}
		\setlength{\tabcolsep}{0pt}
		\begin{tabular*}{\linewidth}{@{\extracolsep{\fill}}cc*{6}{c}@{}}
			\hline
			\multirow{3}{*}{$m$}&\multirow{3}{*}{${1}/{h}$}&\multirow{3}{*}{GMRES}&\multicolumn{5}{c}{$\MatPre$-GMRES}\\\cline{4-8}
			&&&\multirow{2}{*}{Assemble}&\multicolumn{2}{c}{Direct}&\multirow{2}{*}{GMG}&\multirow{2}{*}{BoomerAMG}\\\cline{5-6}
			&&&&Factor&Solve&&\\\hline
			\multirow{4}{*}{$1$}&4&0.02&0.03&--&0.01&--&--\\
			&8&0.31&0.22&0.02&0.06&0.05&0.05\\
			&16&6.36&1.80&1.34&1.03&0.42&0.49\\
			&32&124.48&14.66&261.95&16.78&3.97&4.54\\\hline
			\multirow{4}{*}{$2$}&4&0.04&0.03&--&0.01&0.01&0.01\\
			&8&0.74&0.26&0.02&0.14&0.13&0.14\\
			&16&13.19&1.81&1.34&1.82&1.22&1.27\\
			&32&204.03&14.75&262.47&24.44&10.45&10.85\\\hline
			\multirow{4}{*}{$3$}&4&0.07&0.03&--&0.03&0.03&0.03\\
			&8&1.90&0.23&0.02&0.38&0.36&0.38\\
			&16&40.29&2.97&1.33&4.72&3.76&3.90\\
			&32&668.68&14.78&263.22&55.13&31.98&33.65\\\hline
		\end{tabular*}
		\begin{tablenotes}
			\footnotesize
			\item \makebox[4em][l]{--} indicates that the computational time is below the displayed precision.
		\end{tablenotes}
	\end{threeparttable}
\end{table}

Overall, the numerical results show that, for a fixed approximation degree $m$, the iteration counts of the $\MatPre$-preconditioned GMRES method are relatively small and remain nearly constant as the mesh is refined, confirming the mesh robustness of the proposed preconditioner. Furthermore, employing GMG and BoomerAMG as approximate block solvers results in relatively short CPU solution times, particularly for large-scale discrete systems.

\subsection{The Five-Spot Problem}
In the previous subsection, we investigated the convergence behavior of the preconditioned iterative method with respect to the mesh size $h$. In this subsection, we set the polynomial degree to $\mv=\mq=2$ and further consider a five-spot problem, focusing on the iterative behavior for different meshes and penalty parameters. We also test the proposed method on a nonconforming mesh with hanging nodes.

\begin{example}
	\label{ex:FiveSpot}
	In this example, the computational domain is $\od=(0,1.854075)\times(0,1.854075)$. The sink and source wells are located in the corner elements containing $(0,0)$ and $(1.854075,1.854075)$, denoted by $K_{\mathrm{sink}}$ and $K_{\mathrm{source}}$, respectively. The source/sink term is defined by
	\begin{align*}
		f&= \begin{cases}-\frac{1}{4|K_{\mathrm{sink}}|},&\mathrm{in}\;K_{\mathrm{sink}},\\\frac{1}{4|K_{\mathrm{source}}|},&\mathrm{in}\;K_{\mathrm{source}},\\0,&\mathrm{otherwise},\end{cases}
	\end{align*}
	where $|K_{\mathrm{sink}}|$ and $|K_{\mathrm{source}}|$ denote the areas of the corresponding well elements. The homogeneous boundary condition is given by $\boldsymbol{u}\cdot\boldsymbol{n}=0$ on $\gd$.
\end{example}

\begin{figure}[!t]
	\centering
	\begin{subfigure}[t]{0.45\textwidth}
		\centering
		\includegraphics[width=0.85\textwidth]{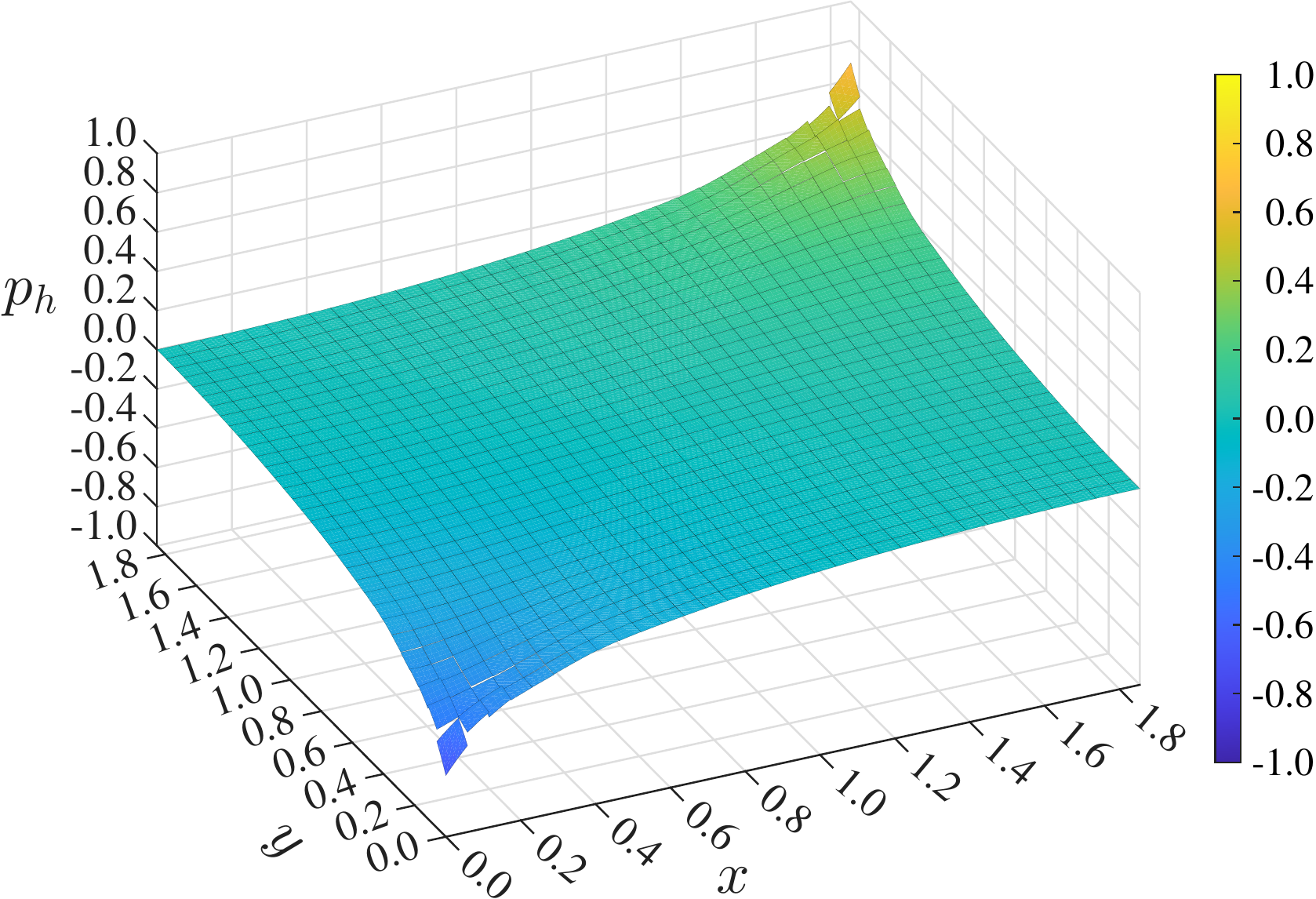}
		\caption{Uniform mesh}
	\end{subfigure}
	\begin{subfigure}[t]{0.45\textwidth}
		\centering
		\includegraphics[width=0.85\textwidth]{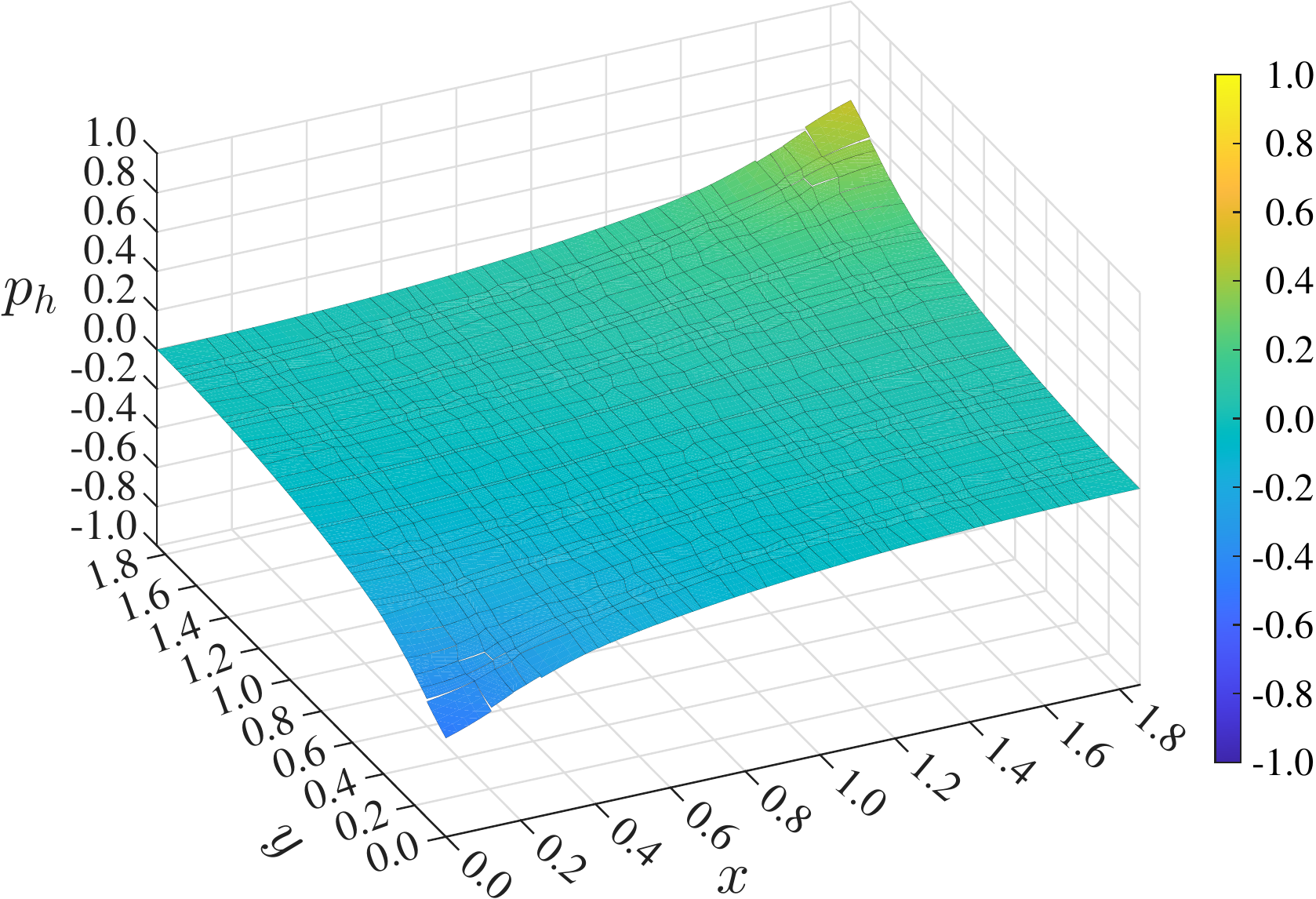}
		\caption{Nonuniform mesh}
	\end{subfigure}
	\begin{subfigure}[t]{0.45\textwidth}
		\centering
		\includegraphics[width=0.85\textwidth]{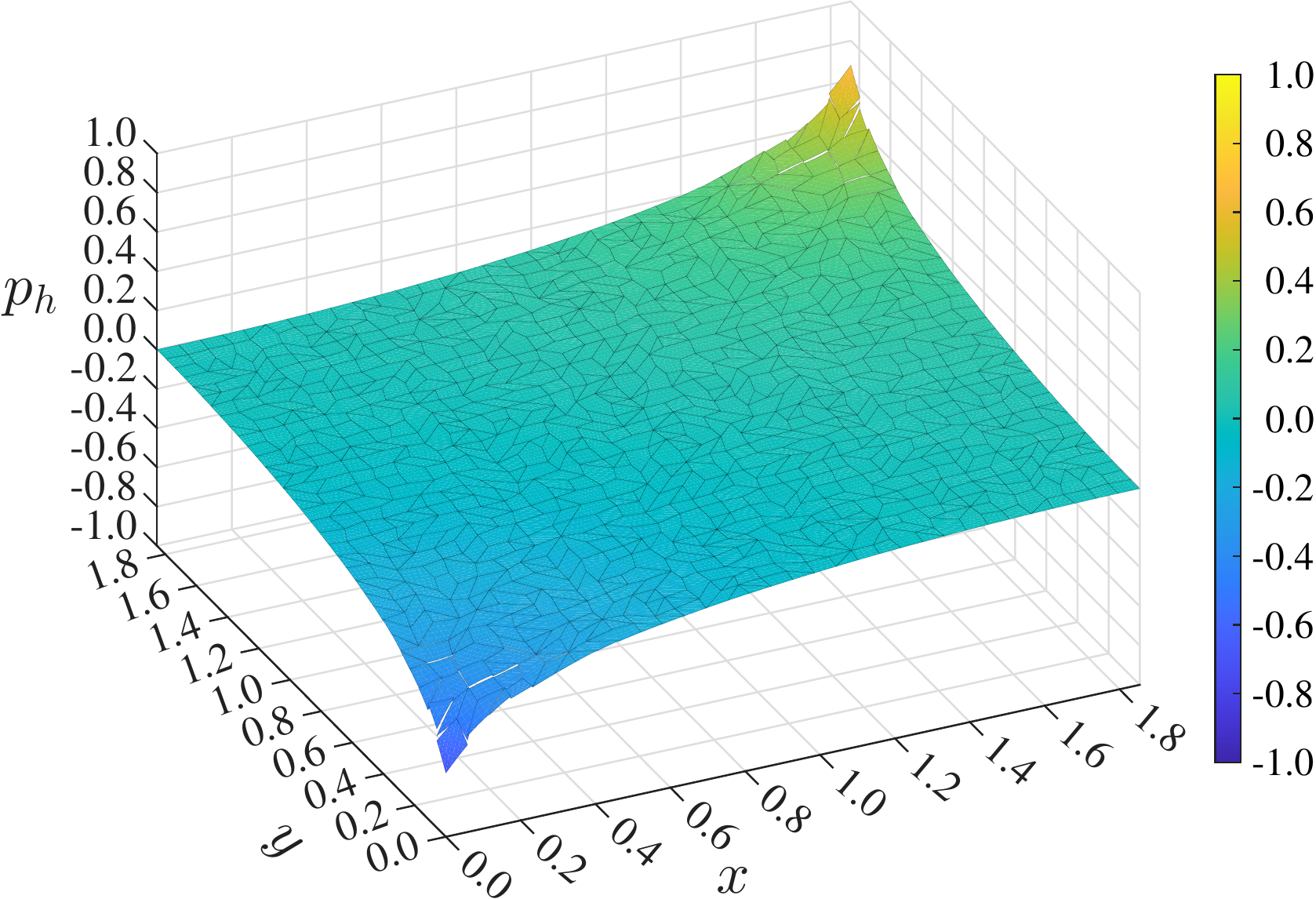}
		\caption{Mixed mesh}
	\end{subfigure}
	\begin{subfigure}[t]{0.45\textwidth}
		\centering
		\includegraphics[width=0.85\textwidth]{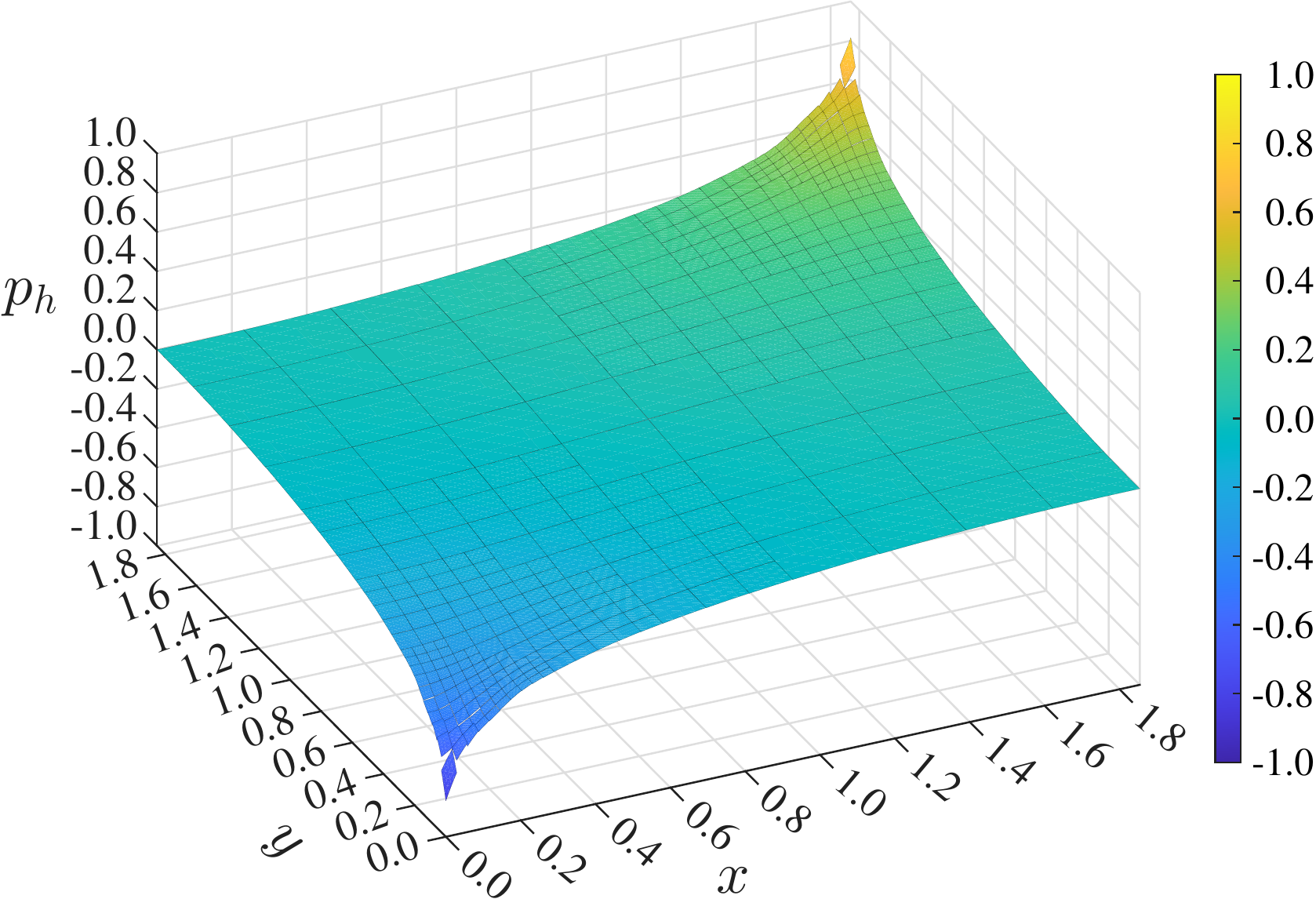}
		\caption{Nonconforming mesh}
	\end{subfigure}
	\caption{Numerical pressure solutions on different meshes in Example~\ref{ex:FiveSpot}.}
	\label{fig:FiveSpot}
\end{figure}

To investigate the influence of mesh geometry, we consider four different meshes: a uniform quadrilateral mesh, a nonuniform quadrilateral mesh, a mixed polygonal mesh, and a nonconforming mesh. The corresponding mesh partitions together with the numerical pressure solutions are presented in Figure~\ref{fig:FiveSpot}.

\begin{figure}[!t]
	\centering
	\begin{subfigure}[t]{0.475\textwidth}
		\centering
		\includegraphics[width=0.65\textwidth]{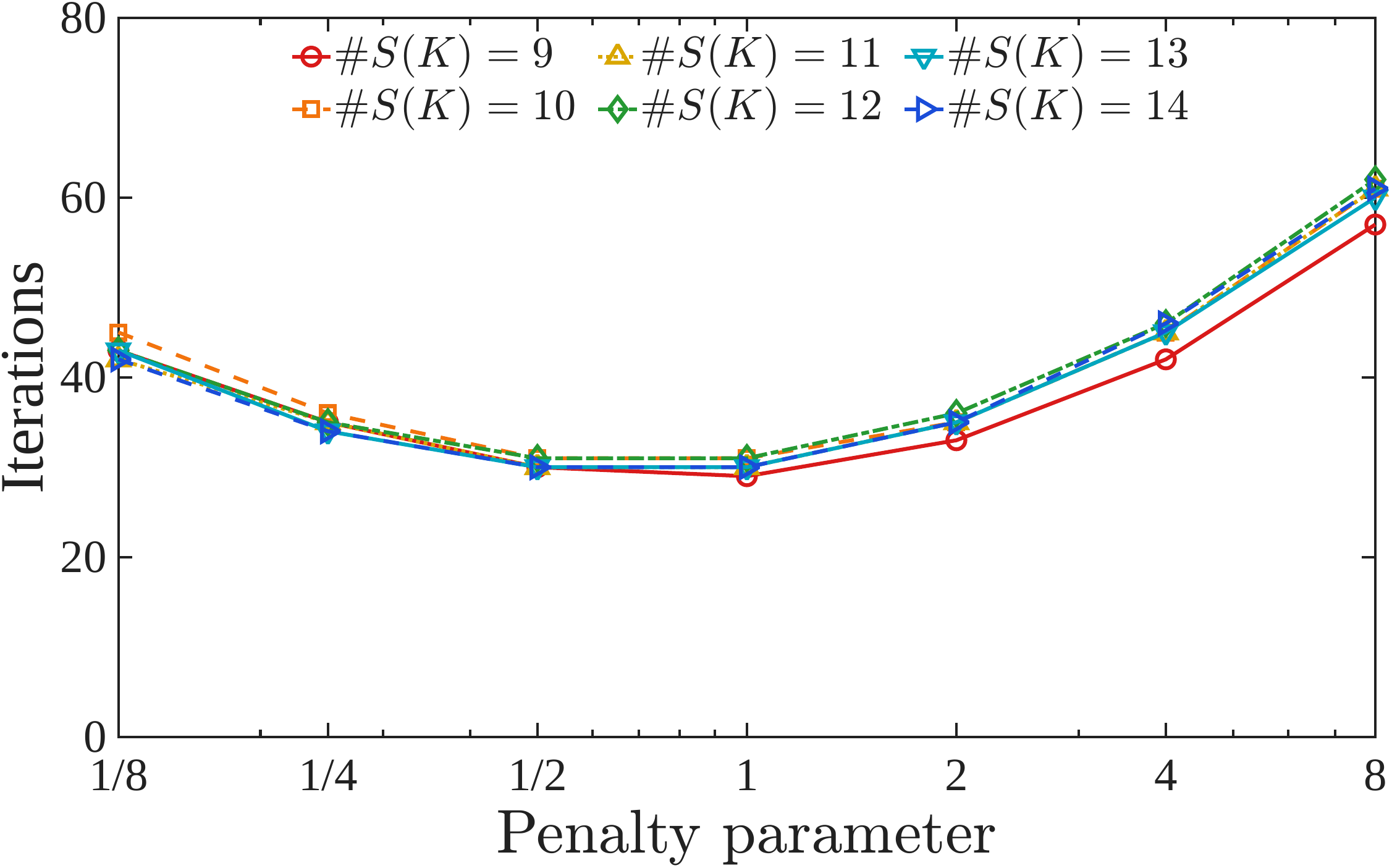}
		\caption{Uniform mesh}
	\end{subfigure}
	\begin{subfigure}[t]{0.475\textwidth}
		\centering
		\includegraphics[width=0.65\textwidth]{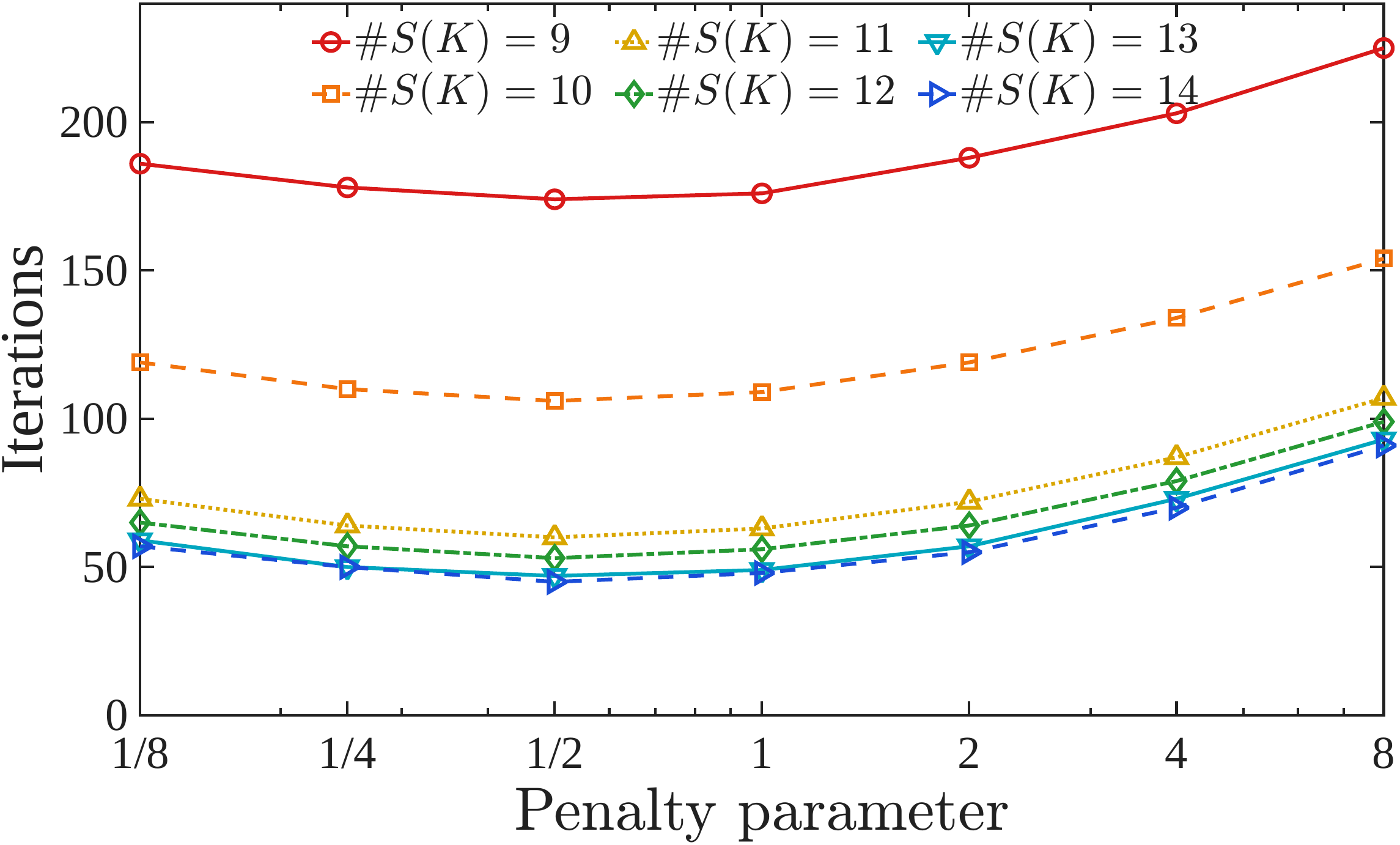}
		\caption{Nonuniform mesh}
	\end{subfigure}
	\begin{subfigure}[t]{0.475\textwidth}
		\centering
		\includegraphics[width=0.65\textwidth]{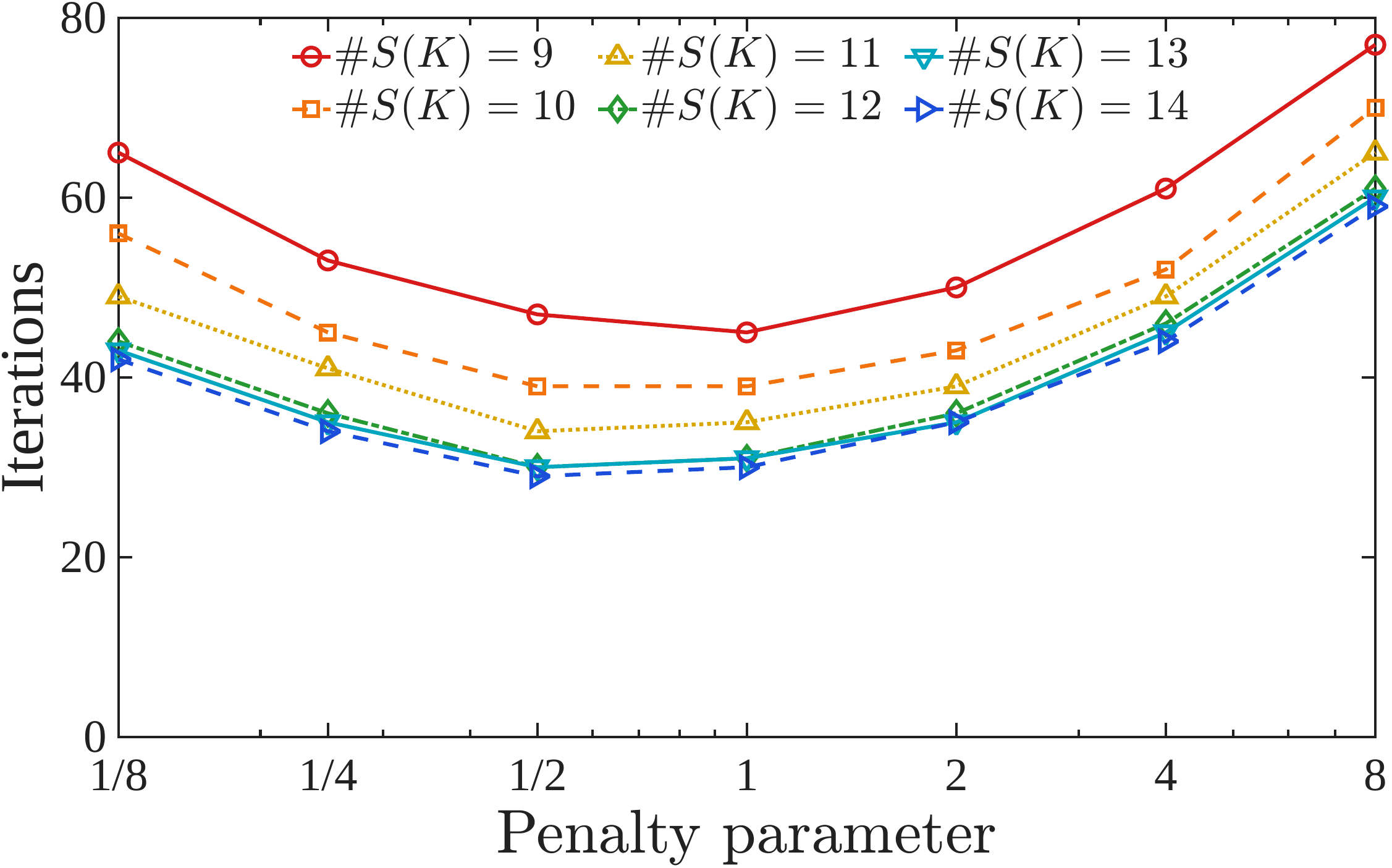}
		\caption{Mixed mesh}
	\end{subfigure}
	\begin{subfigure}[t]{0.475\textwidth}
		\centering
		\includegraphics[width=0.65\textwidth]{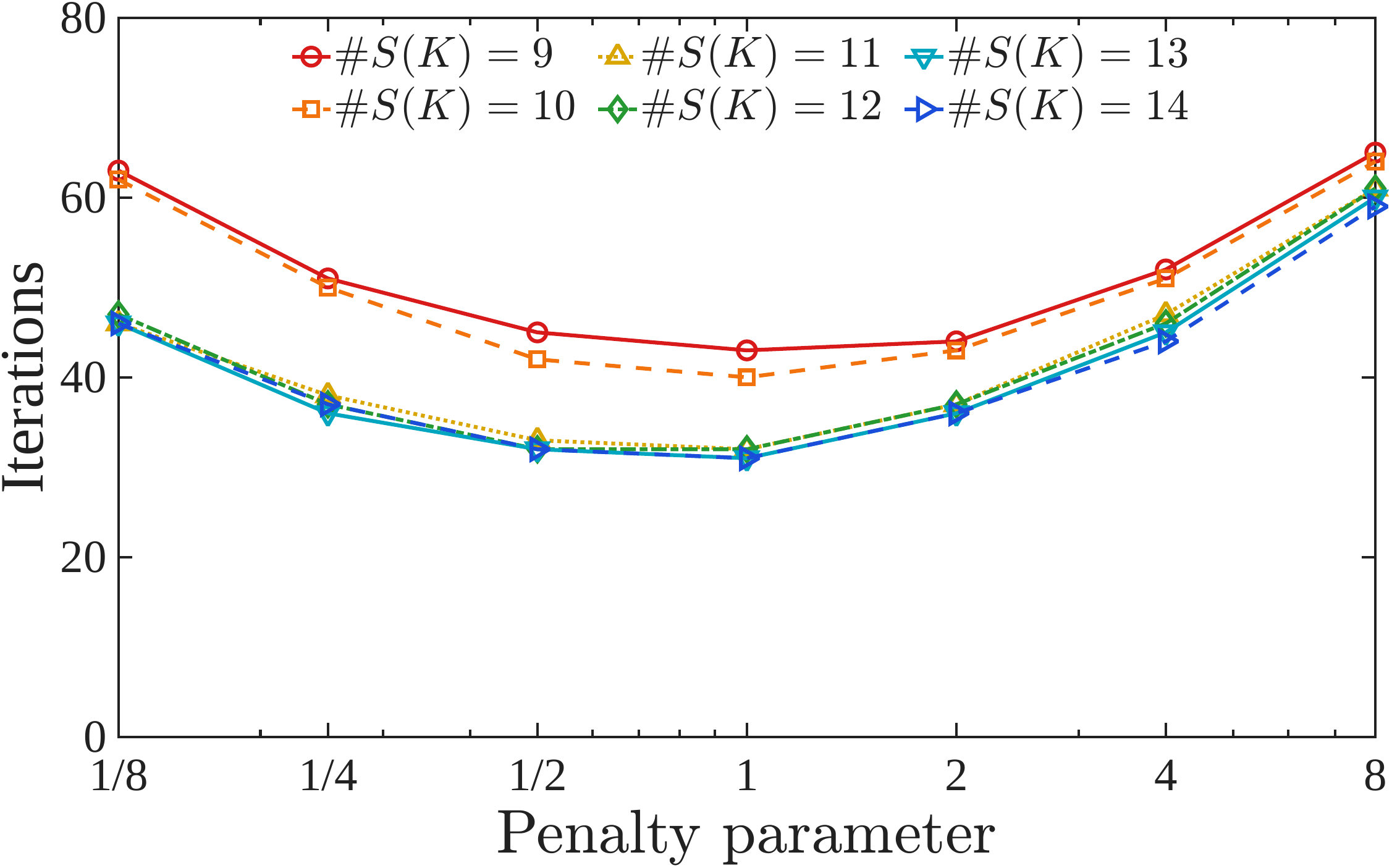}
		\caption{Nonconforming mesh}
	\end{subfigure}
	\caption{GMRES iteration counts for different meshes in Example~\ref{ex:FiveSpot}.}
	\label{fig:FiveSpot-iterations}
\end{figure}

Despite the differences in mesh structures, all the meshes effectively capture the main features of the pressure field. We next set $\alpha=\beta$ and investigate the effect of the penalty parameter on the iterative performance.

Figure~\ref{fig:FiveSpot-iterations} shows the influence of the element patch size and the penalty parameter on the iteration counts for different meshes. The iteration counts exhibit similar trends across the different meshes as the penalty parameter varies. Within the tested range, the iteration counts decrease with increasing penalty parameter, reaching a minimum at $1/2$ or $1$, and then rise as the parameter increases further.

For the uniform quadrilateral mesh, the curves corresponding to different patch sizes nearly overlap, indicating that the patch size has little influence on the iteration counts. In contrast, the nonuniform, mixed, and nonconforming meshes show some sensitivity to the patch size: smaller patches lead to larger iteration counts, while a slight increase in patch size can substantially reduce the iteration counts. This may be related to the variation in element sizes, which can affect the selection of local element patches for the reconstruction. However, further increasing the patch size reduces the iteration counts only slightly, while the cost of reconstruction and matrix assembly continues to increase. Therefore, the patch size should be chosen according to the mesh geometry to balance iterative performance and computational cost.

\subsection{Variable Permeability Problem}
The previous tests demonstrate that the preconditioned solver exhibits robust convergence behavior under mesh refinement. In this subsection, we consider the quadratic case, i.e., $\mv=\mq=2$, and fix the computational mesh. We investigate the effect of different permeability coefficients $\kappa$ on the performance of the proposed method.

\begin{example}\label{ex:Permeability}
	Let $\od=(0,1200)\times(0,2200)\times(0,170)$. The exact pressure is given by
	\begin{align*}
		\pd(x,y,z)=\left(\frac{x}{1200}\right)^2\left(1-\frac{x}{1200}\right)^2\left(\frac{y}{2200}\right)^2\left(1-\frac{y}{2200}\right)^2\left(\frac{z}{170}\right)^2\left(1-\frac{z}{170}\right)^2-\frac{1}{27000},
	\end{align*}
	and the corresponding exact velocity field is defined by
	\begin{align*}
		\ud(x,y,z)=-\kappa(x,y,z)\nabla\pd(x,y,z).
	\end{align*}
\end{example}

In this example, we fix a mesh derived from the SPE10 grid~\cite{Christie2001SPE10}, consisting of $1{,}122{,}000$ hexahedral elements arranged in a regular $60\times220\times85$ Cartesian grid, and investigate the effect of different isotropic permeability coefficients on the performance of the preconditioned iterative solver. Specifically, we consider the following three cases:
\begin{enumerate}[label=\textbf{Case \arabic*.},labelindent=\parindent,leftmargin=*]
	\item $\kappa(\x)=\kappa_{\mathrm{const}}$;
	\item $\kappa(\x)=2+\sin(\frac{2\pi{x}}{1200})\sin(\frac{2\pi{y}}{2200})\sin(\frac{2\pi{z}}{170})$;
	\item $\kappa(\x)=\kappa_{\mathrm{const}}^{\frac{1}{3}\left(\frac{x}{1200}+\frac{y}{2200}+\frac{z}{170}\right)}$,
\end{enumerate}
where $\kappa_{\mathrm{const}}=10^{-6},10^{-4},10^{-2},1,10^{2},10^{4},10^{6}$ for Cases 1 and 3.

\begin{table}[!t]
	\caption{GMRES iteration counts with BoomerAMG and CPU times (s) for different permeability coefficients in Example~\ref{ex:Permeability}.}
	\label{tab:permeability-steps-times}
	\renewcommand{\arraystretch}{1.2}
	\centering
	\begin{tabular}{cc*{3}{>{\centering\arraybackslash}p{5em}}}
		\hline
		&$\kappa_{\mathrm{const}}$&Steps&Assemble&Solve\\\hline
		\multirow{7}{*}{Case 1 }
		&$10^{-6}$&62 &150.27 & 101.68\\
		&$10^{-4}$&62 &156.33 & 102.21\\
		&$10^{-2}$&63 &151.92 & 103.49\\
		&$1$ &63 &150.77 & 103.33\\
		&$10^{2}$ &76 &152.97 & 125.54\\
		&$10^{4}$ &96 &151.18 & 159.62\\
		&$10^{6}$ &97 &150.80 & 160.18\\\hline
		\multirow{1}{*}{Case 2 }
		& -- & 64 &156.03 & 104.30 \\\hline
		\multirow{6}{*}{Case 3 }
		&$10^{-6}$& 63 &154.24 & 103.51\\
		&$10^{-4}$& 61 &155.44 & 100.34\\
		&$10^{-2}$& 57 &154.32 & 94.61\\
		&$10^{2}$ & 61 &156.45 & 100.44\\
		&$10^{4}$ & 75 &153.56 & 122.18\\
		&$10^{6}$ & 89 &155.25 & 147.86\\\hline
	\end{tabular}
\end{table}

Table~\ref{tab:permeability-steps-times} presents the GMRES iteration counts and CPU times for the three permeability cases. The results with $\kappa_{\mathrm{const}}=1$ in Case~3 are omitted from the table, as this case is identical to Case~1 with $\kappa_{\mathrm{const}}=1$. For Case~1, the iteration counts and GMRES solve times remain relatively stable as the constant permeability parameter varies from $10^{-6}$ to $1$, while both increase at larger permeability values. Assembly of the low-order preconditioning matrix requires nearly constant CPU time, showing little dependence on $\kappa_{\mathrm{const}}$ across the tested parameter range.

For Case~2, the permeability coefficient varies within a limited range, and the iteration count and GMRES solve time remain comparable to those for Case~1 with $\kappa_{\mathrm{const}}=1$. The preconditioning matrix assembly time is also comparable. These results indicate that the spatial variation of the permeability coefficient in Case~2 does not substantially increase the computational cost of the preconditioned solve. Despite the exponential variation of the permeability coefficient, Case~3 exhibits iteration counts and GMRES solve times comparable to those of Case~1 across the tested values of $\kappa_{\mathrm{const}}$.

Overall, the GMRES iteration counts depend on the permeability parameter but remain moderate across all tested cases. Notably, even for the SPE10 mesh, the $\MatPre$-preconditioned GMRES method with BoomerAMG approximate block solvers requires at most about 160 s of GMRES solve time. These results demonstrate the effectiveness and computational efficiency of the proposed preconditioned solver for both constant and variable permeability coefficients.

	% vim:spell:tw=70:fo+=Mn:cc=70
\section{Conclusion}\label{sec:Conclusion}

In this paper, we adopt the RDA framework for a stabilized DG discretization of Darcy flow and develop an efficient preconditioning method based on the piecewise constant space. The theoretical analysis derives an optimal-order error estimate in the energy norm and a suboptimal error estimate in the $L^2$ norm, and provides mesh-independent field-of-values bounds for the preconditioned matrix. Finally, two- and three-dimensional numerical experiments verify the convergence of the discretization and demonstrate the effectiveness and computational efficiency of the proposed preconditioning method.

  \section*{Acknowledgements}
	Ming Cui's work was supported by the National Natural Science Foundation of China (Grant No.~12371386). Fanyi Yang's work was supported by the Natural Science Foundation of Sichuan Province (Grant No.~2026NSFSC0023).
	\bibliographystyle{amsplain-nodash}
	\bibliography{References.bib}
\end{document}